\documentclass[graybox]{svmult}

\usepackage{type1cm}        
\usepackage{makeidx}         
\usepackage{graphicx}        
\usepackage{multicol}        
\usepackage[bottom]{footmisc}

\usepackage{newtxtext}       %
\usepackage[varvw]{newtxmath}       

\makeindex             

\let\Bbb\mathbb
\let\goth\mathfrak

\def\og{\leavevmode\raise.3ex\hbox{$\scriptscriptstyle\langle\!\langle\,$}}

\def\fgf{\/\leavevmode\raise.3ex\hbox{$\scriptscriptstyle\,\rangle\!\rangle$}}

\def\fg{\fgf\ }

\newcommand{\carre}{\qed}

\newcommand{\findem}{\ensuremath\blacksquare}

\let\stz\ss

\def\Demd#1|{\parindent=0pt\par{\sl D\'emonstration d#1}\pointir\parindent=20pt}

\let\petcap\sc

\def\Defnns{\parindent=0pt\par{\rm\bf D\'efinitions et notations.\ }\parindent=20pt}

\def\finc{\vskip12pt}

\def\Defn{\parindent=0pt\par{\rm\bf D\'efinition.\ }\parindent=20pt}

\def\Dem{\parindent=0pt\par{\sl D\'emonstration}\pointir\parindent=20pt}

\def\Demo#1|{\parindent=0pt\par{\sl D\'emonstration #1}\pointir\parindent=20pt}

\def\Defns{\parindent=0pt\par{\rm\bf D\'efinitions.\ }\parindent=20pt}

\makeatletter
\newdimen\indentTh\indentTh=0pt
\def\p@int{{\rm .}}
\def\p@intir{\discretionary{\rm .}{}{\rm .\kern.35em---\kern.7em}}
\def\pointir{\afterassignment\pointir@\let\next=}
\def\pointir@{\ifx\next\par\p@int\else\p@intir\fi\next}
\long\def\Thc#1|#2\finc{\Th{}{#1}{\pointir}{#2}}
\long\def\Th#1#2#3#4{\parindent=\indentTh\par\vskip5pt
{#1}{\petcap #2}{\sl #3}\parindent=20pt{\sl #4\par}\vskip 5pt\parindent=20pt}

\newdimen\indentssec\indentssec=20pt
\newdimen\indentrem\indentssec=0pt

\long\def\Remarque#1#2#3#4{\parindent=\indentrem\par\vskip5pt
{#1}{ \sl #2}{\sl #3}\parindent=20pt{#4\par}
\vskip 5pt\parindent=20pt}
\long\def\Rmc#1|#2\finc{\Remarque{}{#1}{\pointir}{#2}}

\long\def\rmc#1|#2\finc{\remarque{}{#1}{\pointir}{#2}}
\long\def\parc#1\finc{\remarque{}{}{}{#1}}
\long\def\remarque#1#2#3#4{\parindent=\indentrem\par\vskip5pt
{#1}{\sl #2}{\sl #3}\parindent=20pt{#4\par}
\vskip 5pt\parindent=20pt}
\makeatother

\def\Demdsp#1|{\parindent=0pt\par{\sl D\'emonstration d#1.}\parindent=20pt}
\def\oldstyle{}

\begin{document}

 
\title*{$\displaystyle SL(2, {\Bbb Z})$%
, les tresses \`a trois brins, le tore modulaire et
$Aut^{+}(F_{2})$%
.}

\author{{\it par} Alexis Marin {\it et d'apr\`es} Emil Artin {\it et\/} Jacob Nielsen}
\maketitle
\vskip-30mm
\hfill { \`a Valentin Po\'enaru pour son nonante et deuxi\`eme anniversaire,\break
  \null\hfill de la part de son ami Jean Cerf%
  \footnote{Ce texte a d'abord \'et\'e accept\'e dans un volume \'edit\'e par Louis Funar et Athanase Papadopoulos en l'honneur de Valentin Poenaru, puis refus\'e par Springer pour ne pas \^etre \'ecrit dans la \og grande belle langue angloUSienne\fgf.}.\/}

\vskip15mm


\centerline{ R\'esum\'e}
\vskip3mm
    {
      \sl
Le groupe sp\'ecial lin\'eaire entier deux-dimensionnel agit 
 sur le quotient par sym\'etrie centrale du tore entier.
La pr\'esentation d'Artin
du groupe
$B_{3}$
des tresses \`a trois brins produit alors
une pr\'esentation de
$SL(2, {\bf Z})$
de g\'en\'erateurs, les paraboliques
$\begin{pmatrix}1& -1\cr 0& 1\cr\end{pmatrix}$
and
$\begin{pmatrix}1& 0\cr 1& 1\cr\end{pmatrix}$%
. Cette
\og pr\'esentation d'Artin\fg
de
$SL(2, {\bf Z})$
clarifie l'action de son groupe d\'eriv\'e
sur le demi-plan de Poincar\'e
${\bf H}$
et son quotient le tore modulaire,
ainsi que  le th\'eor\`eme de Nielsen donnant le groupe des automorphismes directs
du groupe libre de rang
$2$
en produit semi-direct, amalgam\'e
sur le sous-groupe d'indice
$2$
du centre de
$B_{3}$%
, des automorphismes int\'erieurs par
$B_{3}$%
. 
}
\vskip10mm

\abstract{The action of 
$SL(2, {\bf Z})$
on the integer torus and its quotient by central symmetry
and Artin's presentation of three strings braid group
$B_{3}$%
, produces a presentation with parabolic generators
$\begin{pmatrix}1& -1\cr 0& 1\cr\end{pmatrix}$
and
$\begin{pmatrix}1& 0\cr 1& 1\cr\end{pmatrix}$%
. This braided presentation describes the action of the derived group
on Poincar\'e's half plane and its quotient the modular torus, just as Nielsen's theorem
giving the group of  direct automorphisms of the free group on two generators  as
 semi-direct product, amalgamated on the index
 $2$
 subgroup of the center of
 $B_{3}$%
 , of  inner automorphisms with
$B_{3}$%
.
}

\vskip12pt

\noindent 2020 {\sl Mathematics Subject Classification}~:\hfill\break
Primary: 20F36; Secondary: 20H05, 20F05, 20H05, 20E36, 57K20, 01A99.

\def\keywordname{{\bf Key Words and Phrases:}}%

\vskip15mm
Soit
$Aut(F_{2})$
le groupe d'auto\-morphismes de
$F_{2}\!=\!\langle u, v\rangle$%
, le groupe libre \`a deux g\'en\'erateurs, et
$Int(F_{2})$
le sous-groupe normal des automorphismes int\'erieurs.

Le groupe d'Artin des tresses \`a trois brins
$B_{3} = \langle a, b\ |\ aba = bab =: s\rangle$%
, de centre%
\footnote{\small
On n'utilisera que l'\'el\'ement
$s\!=\!aba\!=\!bab$
 conjugue
$a$
en
$b$
et
$b$
en
$a$%
, donc
$s^{2}$
est central.\hfill\break
Que
$s^{2}$
engendre le centre demande plus de travail dans la litt\'erature
mais suivra ici de la pr\'esentation d'Artin du groupe de centre trivial
$PSL(2, {\Bbb Z\/})\!=\!B_{3}/\langle s^{2}\rangle$
\'etablie au \S1.
}
$Z = \langle s^{2}= ababab\rangle$
se repr\'esente aussi%
\footnote{via
$a(u, v) = (u, u^{-1}v),\quad b(u, v) = (vu, v)$
(Cf. la Proposition de l'Appendice B).
}\
comme sous-groupe de
$Aut(F_{2})$%
.

L'ab\'elianisation
$ F_{2}\rightarrow {\Bbb Z}^{2}$
induit
$\rho : Aut(F_{2}) \rightarrow GL(2, {\Bbb Z})$%
. Le  sous-groupe des automorphismes {\it directs\/} de
$F_{2}$
est not\'e~:
$$Aut^{+}(F_{2}) := \rho^{-1}(SL(2, {\Bbb Z}))\subset Aut(F_{2})\ .$$

Nielsen a essentiellement montr\'e~:

\Thc Th\'eor\`eme|
Le sous-groupe
$Int(F_{2})$
des int\'eriomorphismes, 
intersecte
$B_{3}$
en le sous-groupe d'indice deux de son centre
$Int(F_{2})\cap B_{3} = \langle s^{4}\rangle$%
, est le noyau  de~:
$$\rho^{+} := \rho_{|}\!: Aut^{+}\!(F_{2}) \rightarrow SL(2, {\Bbb Z})$$
et engendre avec
$B_{3}$
le groupe des automorphismes directs~:
$$Aut^{+}(F_{2}) = Int(F_{2}) B_{3}\ .$$

\finc

Comme
$\rho^{+}$
est surjectif et
$\rho^{+}(s^{2})\!=\!-I\!\ne\!1$
il suit la {\it pr\'esentation\/} d'Artin~:
$$SL(2, {\Bbb Z}) \simeq B_{3}/\!\langle s^{4}\rangle\, =\,
\langle a, b\ |\ aba = bab, (ababab)^{2} = 1\rangle \leqno{({\goth A})}$$

Mais la preuve propos\'ee en Appendice B de ce th\'eor\`eme de Nielsen
n\'eces\-site cette pr\'esentation,
pr\'esentation qui, ind\'ependamment de la surjectivit\'e de
$\rho^{+}$%
, est obtenue au \S1 gr\^ace \`a l'action affine de
$SL(2, {\Bbb Z})$
sur la sph\`ere plate
${\Bbb S}$%
, quotient du tore entier
${\Bbb T} =  {\Bbb R}^{2}/{\Bbb Z}^{2}$
par la sym\'etrie centrale
$-I: {\Bbb T}\rightarrow {\Bbb T}, x\mapsto -x$%
.

Le \S2
d\'ecrit, \`a l'aide de
$({\goth A})$%
, l'action sur le demi-plan de Poincar\'e du groupe d\'eriv\'e
$G'$
du groupe modulaire
$G = SL(2, {\Bbb Z})/\{\pm I\}$%
, de quotient le {\it tore modulaire.\/}

L'appendice A explicite l'image dans
$B_{3}\!=\!\pi_{0}(Diff({\Bbb S}; {\Bbb T}_{2}, \{\overline{0}\}))$
des g\'en\'e\-rateurs de la pr\'esentation
$({\goth A})$%
, l'appendice B \'etablit, en sus  du plongement de
$B_{3}$
dans
$Aut^{+}(F_{2})$%
,  les th\'eor\`emes de Nielsen de structure et caract\'erisation de
$Aut^{+}(F_{2})$%
, l'appendice C, gr\^ace \`a une forme normale dans
$B_{3}$%
, d\'eduit du th\'eor\`eme de Nielsen la classification \`a conjugaison pr\`es
des \'el\'ements de torsion de
$SL(2, {\Bbb Z\/})$
et
$Aut^{+}(F_{2})$%
, l'appendice D introduit le groupe des {\it tresses di\'edrales\/}
\`a trois brins
$DB_{3}$
permet\-tant d'obtenir des r\'esultats analogues pour le groupe complet
$Aut(F_{2})$
des automorphismes du groupe libre \`a deux g\'en\'erateurs,
l'appendice E traite
du produit semi-direct
de deux groupes amalgam\'e%
\footnote{une notion, sans doute bien connue,
mais que nous n'avons pas trouv\'ee dans la litt\'erature.
}\
sur un sous-groupe commun.

Enfin en postroduction, des commentaires bibliographiques sur les divers
traitements
de ces r\'esultats dans la litt\'erature.

\section*{{\sc Remerciements}}

Cette note doit tout aux relectures attentives de Danielle Bozonat, Daniel Marin, des \'editeurs
Athanase Papdopoulos et Louis Funar
et surtout Greg McShane qui,
en partageant ses tentatives contre
la conjecture d'unicit\'e de Frobenius pour les nombres de Markoff
(Cf. {\bf [H]}),
a pouss\'e le metteur en sc\`ene \`a comprendre le tore modulaire.
Enfin n'oublions pas l'arbitre qui, avec  patience et longueur de temps,
a fourni la mati\`ere du  quatri\`eme de couverture.


\section{Le morphisme $SL(2, {\Bbb Z})\rightarrow B_3/Z$ et la pr\'esentation parabolique de $SL(2, {\Bbb Z})$.}

Le groupe
 $SL(2, {\Bbb Z})$
 agit sur l'espace vectoriel
 $V\!=\!{\Bbb R}^{2}$%
 , respectant le r\'eseau
 $\Lambda\!=\!{\Bbb Z}^{2}$%
 , induisant un automorphisme du tore quotient
$\displaystyle {\Bbb T}\!=\!V/\Lambda$%
.
La suite exacte~:

$$0\rightarrow \Lambda\rightarrow V\rightarrow {\Bbb T}\rightarrow 0$$
fait appara\^{\i}tre
$\Lambda$
comme groupe fondamental de
${\Bbb T}$
et 
$\overline{M}: ({\Bbb T}, 0)\rightarrow ({\Bbb T}, 0)$%
, l'automorphisme induit par
$M\in SL(2, {\Bbb Z})$%
, produit sur ce groupe fondamental la restriction
$\pi_1(\overline{M}, 0) = M_{|} : \Lambda\rightarrow\Lambda$
de
$M$
\`a
$\Lambda$%
.

Le centre
$\{\pm I\} < SL(2, {\Bbb Z})$
fixe le sous-groupe
de
$2$%
-torsion
${\Bbb T}_2 = {{1}\over{2}}{\Bbb Z}/{\Bbb Z}$
et l'action de
$SL(2, {\Bbb Z})$
passe au quotient en une action par hom\'eomorphismes affines~:

$$\pi : G=PSL(2, {\Bbb Z}) \rightarrow
Affeo({\Bbb T}/\{\pm I\}, {\Bbb T}_2,  \{\overline{0}\})$$
du groupe modulaire
$G\!=\!SL(2, {\Bbb Z})/\{\pm I\}$
sur la sph\`ere plate \`a
$3 + 1$
 points singuliers (d'angle
$\pi$%
)%
\footnote{
 not\'es
$p = (\overline{{{1}\over{2}}, 0)},\, n = \overline{(0, {{1}\over{2}})},\,
m = \overline{({{1}\over{2}}, {{1}\over{2}})}$
et
$\infty = \overline{(0, 0)}$%
.
}\
$({\Bbb T}/\{\pm I\}, {\Bbb T}_2, \{0\})
=:({\Bbb S}, \{m, n, p, \infty\}, \{\infty\})=:
({\Bbb S}, R, \{\infty\})$%
.

D'o\`u par lissage%
\footnote{
Si
$M=\begin{pmatrix}a& b\cr c& d\cr\end{pmatrix}$
on d\'eforme, dans le commutant de
$\{\pm I\}$%
,
$\overline{M}$
 au voisinage de
$R$
en
$\overline{M}_c$%
, affine pr\`es de
$R$
de partie lin\'eaire la \og conformis\'ee
$\begin{pmatrix}a& -c\cr c& a\cr\end{pmatrix}$
 sur la premi\`ere colonne\fgf
, 
de
$M$%
.
}%
, fibration de Cerf
et suite exacte  d'une fibration~:

$$\pi_0 : SL(2, {\Bbb Z})=\pi_{0}(SL(2, {\Bbb Z}))\rightarrow
\pi_0(Affeo({\Bbb S}, R,  \{\infty\}))\rightarrow$$
$$\rightarrow\pi_0(Diff({\Bbb S}, R,  \{\infty\}))
= \pi_1(Pl(R\setminus\{\infty\}, {\Bbb S}\setminus\{\infty\}))/Z$$
o\`u
$Z = \langle ababab\rangle$
est le centre du groupe
$B_3$
des tresses \`a trois brins d'Artin%
\footnote{
$a$
(resp.
$b$%
)
 fixe
$p$
(resp.
$n$%
)
et  \'echange, par demi-tour positif,
$n$
(resp.
$p$%
)
et
$m$%
.
}~:

$$B_3 = \pi_1(Pl(\{m, n, p,\}, {\Bbb S}\setminus\{\infty\})) =\
\langle a, b\ |\ aba = bab\rangle\ .$$

Si
$A= \begin{pmatrix}1& -1\cr 0& 1\cr\end{pmatrix}, B = \begin{pmatrix}1& 0\cr 1& 1\cr\end{pmatrix}$
et
$S = \begin{pmatrix}0& -1\cr 1& 0\cr\end{pmatrix}$%
, on a~:

$$ABA\!=\!S\!=\!BAB,\ S^{4}\!=\!I,\ S^{2}\!=\!-I\!\ne\!I\quad
\hbox{\rm et {(cf. Appendice A)}\/}\quad
\pi_0(A)\!=\!a, \pi_0(B)\!=\!b$$
d'o\`u un morphisme
$\sigma : B_3\rightarrow SL(2, {\Bbb Z})$
d\'efini sur les g\'en\'erateurs de la pr\'esentation d'Artin par
$\sigma(a) = A, \sigma(b) = B$
qui, si on note
$\overline{\rho} : B_3\rightarrow B_3/Z$%
, v\'erifie
$\pi_0\circ \sigma = \overline{\rho}$%
.
\vskip2mm

Ainsi, si
$s = aba=bab$%
, le morphisme
$\sigma$
induit un morphisme injectif~:
$$\displaystyle\overline{\sigma} : B_3/\langle s^{4}\rangle\ \rightarrow SL(2, {\Bbb Z})$$
qui est  isomorphisme puisqu'aussi surjectif car,  si 
$M\!\in\!SL(2, {\Bbb Z})$
et
$\pi_0(M)$
est repr\'e\-sent\'e par
$\mu\!\in\!B_3$%
, les matrices 
$M$
et
$\sigma(\mu)$
sont la
$\pi_1$%
-action de relev\'es de deux diff\'eomorphismes
homotopes de
$({\Bbb S}, R, \{\infty\}) = ({\Bbb T}/\{\pm I\}, {\Bbb T}_2, \{0\})$
donc sont, soit \'egales
et
$M = \sigma(\mu)$%
, soit oppos\'ees et
$M = \sigma(s^{2}\mu)$%
.\hfill\carre
\vskip2mm

On vient d'\'etablir les pr\'esentations~:
$$SL(2, {\Bbb Z}) = \langle A, B\ |\ ABA(BAB)^{-1}, (ABABAB)^{2}\rangle$$
$$G = PSL(2, {\Bbb Z}) = \langle A, B\ |\ ABA(BAB)^{-1}, ABABAB\rangle$$
et
$SL(2, {\Bbb Z})$
et
$PSL(2, {\Bbb Z})$
ont leurs ab\'elianis\'es cycliques d'ordre 
$12$
et $6$%
.\hfill\findem


\section{Action du groupe d\'eriv\'e $SL(2, {\Bbb Z})'$ sur le demi-plan de Poincar\'e.}

Soit dans le demi-plan de Poincar\'e
${\Bbb H}=\{z\in{\Bbb C} | \Im(z)>0\}$
les translat\'es
par les puissances de
$A$
du domaine fondamental usuel de l'action de
$G = PSL(2, {\Bbb Z})$%
~:
$${\cal B}_{n} = A^{-n}\bigl(\{z\in{\Bbb C}\ |\
\Im(z)>0, -{{1}\over{2}}\leq \Re(z)\leq {{1}\over{2}}, |z|\geq 1\}\bigr)\ .$$

L'union
$\tilde{H}\!=\!\cup {\cal B}_n$
quotient\'ee par l'image dans
$G$
du sous-groupe
$\langle A^{-6}S^{2}\rangle$
du groupe d\'eriv\'e%
\footnote{
car,  dans l'ab\'elialis\'e on a~:
$\overline{A}=\overline{B}$
donc
$\overline{S}=\overline{A}^{3}$%
, ainsi
$\overline{A^{-6}S^{2}}=\overline{1}$%
.}\
$SL(2, {\Bbb Z})'$
 de
$SL(2, {\Bbb Z})$
est un hexagone hyperbolique point\'e qui,
par identification de ses c\^ot\'es oppos\'es,
produit un tore sym\'etrique point\'e
${\cal T}$%
, identifications par les images dans
$G$
des \'el\'ements suivants du groupe d\'eriv\'e~:

$$f_n = A^{-n}A^{-3}SA^{n} = A^{-(n + 3)}SA^{n}
= A^{-(n + 1)}A^{-1}BA^{(n + 1)}$$
$$f_{n + 1}f_{n - 1} = A^{-(n + 4)}SA^{(n + 1)}\!A^{-(n + 2)}\!SA^{n-1}
= A^{-(n + 3)}BAA^{(n + 1)}\!A^{-(n + 2)}SA^{n-1}$$
$$= A^{-(n + 3)}BSA^{n-1} = A^{-(n + 3)}SAA^{n-1} = A^{-(n + 3)}SA^{n} = f_n\ .$$
De plus 
$A^{-6}S^{2}\!=\!f_{n} f_{n-3}\!=\!(f_{n-1}f_{
n-2}^{-1})(f_{n-1}^{-1}f_{n-2})\!=\!%
[f_{n-1}^{-1}, f_{n-2}]$%
.

Ainsi%
\footnote{
notant par abus
$f_{n}\in G$
les images dans le groupe modulaire
$G$
des
$f_{n}\in SL(2, {\Bbb Z})$%
.}\
${\cal T}$
est quotient de
${\Bbb H}$
par le sous-groupe de
$G$
engendr\'e par les
$f_{n}$%
~:

$$\langle f_n; n\in{\Bbb Z}\rangle\, =\, \langle f_{-2}, f_{-1}\rangle$$
qui est%
\footnote{
car
$\cup_{n=0}^{5}{\goth B\/}_{n}$%
, union de
$6$
translat\'es du domaine fondamental usuel de
$SL(2, {\Bbb Z\/})$%
, est un domaine fondamental de ce sous-groupe
$\langle f_n; n\in{\Bbb Z}\rangle$%
.
}\
d'indice
$6$
dans le groupe modulaire
$G$
(l'indice de son l'ab\'elianis\'e
$G'$%
).

Comme
$\langle f_{n}; n\in{\Bbb Z}\rangle\,\subset G'$%
, ce sous-groupe est le groupe d\'eriv\'e de
$G$%
~:
$$G'=PSL(2, {\Bbb Z})'\,=\, \langle f_{-2}, f_{-1}\rangle$$
et le tore modulaire est
${\cal T}\!$%
, de parabolique de cusp,
le commutateur%
\footnote{
avec la convention
$[a, b] = a^{-1}b^{-1}ab$
de Bourbaki, Hall\dots
}\
de
$f_{-1}^{-1}$ et
$f_{-2}$%
.\
\null\hfill\carre

Comme
$-I=S^{2}\in SL(2, {\Bbb Z})$
est d'ab\'elianis\'e
$\overline{A}^{6}\ne0,\ -I\not\in SL(2, {\Bbb Z})'$%
, ainsi le quotient
$SL(2, {\Bbb Z})\rightarrow G$
est injectif sur le groupe d\'eriv\'e
$SL(2, {\Bbb Z})'$ 
et ce groupe d\'eriv\'e est engendr\'e%
\footnote{librement puisque c'est le groupe fondamental d'un tore point\'e.
\ }\
par
$f_{-2}\!=\!BA^{-1}\!=\!\begin{pmatrix}1& 1\cr 1& 2\cr\end{pmatrix}$
et
$f_{-1}\!=\!A^{-1}B\!=\!\begin{pmatrix}2& 1\cr 1& 1\cr\end{pmatrix}$%
.
\hfill\findem


\section*{APPENDICE A\\
Description de
${\Bbb T}/\{\pm I\}$
et v\'erification de
$\pi_0(A)=a, \pi_0(B)=b$%
.
}

En identifiant
$(V, \Lambda)$
au plan complexe muni du r\'eseau de Gau\stz
$\,({\Bbb C}, {\goth G})$
o\`u
${\goth G}\!=\!{\Bbb Z}+i{\Bbb Z}$%
, la fonction
$\wp : {\Bbb C}\rightarrow{\Bbb C}$
de Weierstra\stz~:
$$\wp(z)={{1}\over{z^{2}}}+\sum_{\gamma\in{\goth G}\setminus\{0\}}
{{1}\over{(z-\gamma)^{2}}} -{{1}\over{\gamma^{2}}}$$
paire et invariante par
${\goth G}$%
, identifie
${\Bbb T}/\{\pm I\}$
\`a
${\Bbb C}\cup\{\infty\}$
avec
$\wp(\overline{0})=\infty$
et%
\footnote{comme
$\wp(iz)\!=\!-\wp(z)$
et
$m-im\in\Lambda$
on a
$\wp(m)\!=\!0$
et
$\wp(n)\!=\!-\wp(p)$%
. Ainsi
$\wp^{-1}(0)\!=\!m+\Lambda$%
, donc comme
$\wp(\overline{z})\!=\!\overline{\wp(z)}$
et,
pr\`es de
$0$%
, on a
$\wp(z)\sim{{1}\over{z^{2}}}$%
, on a
$\wp(]0, p])\subset ]0, +\infty[$%
.
}~:
$$\wp(m)=\wp({{1}\over{2}}(1+i))=0,\
\wp(p)=\wp({{1}\over{2}})=-\wp({{1}\over{2}}i)=-\wp(n)\in\, ]0, +\infty[\ .$$

Les deux \og sous-groupes de coordonn\'ees\fg
${\bf t}={\Bbb R}/{\Bbb Z},\ {\bf u}={\Bbb R}i/{\Bbb Z}i$
et les deux
\og sous-groupes diagonaux\fg
$\Delta_{+}={\Bbb R}(1+i)/{\Bbb Z}(1+i)$
et
$\Delta_{-}={\Bbb R}(1-i)/{\Bbb Z}(1-i)$%
, points fixes d'isom\'etries conformes%
\footnote{$z\mapsto \pm \overline{z}$
et
$z\mapsto \pm i\overline{z}$
qui toutes deux fixent l'origine
$0$%
.
}\
du r\'eseau
${\goth G}$
commutant \`a
$-I$%
, vont dans des droites et sont d'images%
\footnote{remarquer de nouveau que, pr\`es de
$0$%
, on a
$\wp(z)\sim{{1}\over{z^{2}}}$%
.}\
les demi-axes r\'eels tronqu\'es et  les demi-axes imaginaires~:
$$\wp({\bf t})\!=\!H_{\wp(p)}^{+}\!:=\!%
[\wp(p), +\infty[\cup\{\infty\},\
\wp({\bf u})\!=\!%
H^{\wp(n)}_{-}\!:=\, ]-\infty, \wp(n)]\cup\{\infty\}\subset{\Bbb R}\cup\{\infty\}$$
$$\wp(\Delta_{+})\!=\!%
V_{-}\!: =\, ]-\infty, 0]\,i\cup \{\infty\}, \ 
\wp(\Delta_{-})\!=\!V^{+}\!:=\, [0, +\infty[\,i\cup\{\infty\}
\subset {\Bbb R}i\cup\{\infty\}$$

Ainsi l'action de
$A=\begin{pmatrix}1& -1\cr 0& 1\cr\end{pmatrix}$
fixe
$H^{+}_{\wp(p)}$
et envoie
$H_{-}^{\wp(n)}, V_{-}$
sur 
$V^{+}, H_{-}^{\wp(n)}$
 respectivement, c'est le \og demi-tour\fg positif
$a$%
.

Et l'action de
$B=\begin{pmatrix}1& 0\cr 1& 1\cr\end{pmatrix}$
fixe
$H_{-}^{\wp(n)}$
et envoie
$H^{+}_{\wp(p)}, V^{+}$
sur
$V_{-}, H^{+}_{\wp(p)}$
 respectivement, c'est le \og demi-tour\fg positif
$b$%
.
\hfill\findem


\section*{APPENDICE B\\
Les th\'eor\`emes de Nielsen de structure
et caract\'erisation de
$Aut^{+}(F_{2})$%
.
}

Soit
$F_{2} =\, \langle u, v\rangle$
le groupe libre \`a deux g\'en\'erateurs not\'es
$u, v$%
.

On identifie le groupe libre
$F_{2}$
\`a son groupe d'automorphismes int\'erieurs par~:
$$\varphi: F_{2}\rightarrow Int(F_{2}),\
x\mapsto \varphi_{x}: y\mapsto xyx^{-1}\!=
y^{x^{-1}}\!=:\, ^{x}y$$
et note
$Y\!=\!\langle s^{4}\rangle$
le sous-groupe d'indice
$2$
du centre
$Z\!=\!\langle s^{2}\rangle$
de
$B_{3}$
identifi\'e, par
$s^{4}\mapsto \varphi_{[u, v]}\in Int(F_{2})$%
, au sous-groupe cyclique infini
$Y\!=\!\langle \varphi_{[u, v]}\rangle\subset Int(F_{2})$%
. 

\Thc Proposition| {\sl Les endomorphismes
$\alpha, \beta : F_{2}\rightarrow F_{2}$
d\'efinis par~:
$$\alpha(u, v) = (u, u^{-1}v)\quad
\beta(u, v) = (vu, v)$$
sont des isomorphismes d'inverses d\'efinis par~:
$$\alpha^{-1}(u, v) = (u,  uv),\quad \beta^{-1}(u, v) = (v^{-1}u, v)$$
v\'erifiant la relation de tresse
$\alpha\beta\alpha = \beta\alpha\beta=:s\,:
(u, v)\mapsto (v^{u}, u^{-1})$%
, avec de plus~:
$s^{2}\,: (u, v)\mapsto (u^{-1}, v^{-1})^{vu}$
et
$ s^{4}\,: (u, v)\mapsto (u, v)^{v^{-1}u^{-1}vu}\!=\!(u, v)^{[v, u]}\!=\!%
\varphi_{[u, v]}(u, v)$%
.

Ainsi
$\alpha, \beta$
d\'efinissent un morphisme injectif~:
$$\Psi: B_{3}\rightarrow Aut(F_{2}),\ \Psi(a) = \alpha, \Psi(b) = \beta\ .$$
}
\finc
Par abus, on identifiera
$B_{3}$
au sous-groupe
$\Psi(B_{3})$
de
$Aut(F_{2})$%
.  Soit~:
$$\psi: B_{3}\rightarrow Aut(Int(F_{2})),\
\psi(x): h \mapsto x\circ h\circ x^{-1}\ .$$                
Le sous-groupe 
$H\!=\!Int(F_{2})B_{3}\!\subset\!Aut^{+}(F_{2})$
est%
\footnote{Cf. {\bf Corollaire\/} de l'appendice E.
}\
produit semi-direct   suivant
$\psi$~:
$$H=Int(F_{2})\rtimes_{\psi, Y}B_{3}\subset Aut^{+}(F_{2})$$
de
$Int(F_{2})$
par
$B_{3}$
amalgam\'e sur
$Y$
(identifiant
$\varphi_{[u, v]}$
\`a
$s^{4}$%
).

\vskip2mm
\Demd e la Proposition| Que
$\alpha^{-1}$
soit inverse de
$\alpha$
suit de~: 
$$
\alpha\alpha^{-1}(u, v) = \alpha(u, uv)
= (u, uu^{-1}v) = (u, v)
$$
donc
$\alpha\alpha^{-1} = Id$%
, de m\^eme on \'etablit
$\alpha^{-1}\alpha = Id = \beta\beta^{-1} = \beta^{-1}\beta$%
.
\hfill\carre

Pour les relations, notant
$s=\alpha\beta\alpha,\, s'=\beta\alpha\beta$%
, on a~:
$$s(u, v) = \alpha\beta(u, u^{-1}v)\!=\!%
\alpha(vu, u^{-1}v^{-1}v)\!=\!\alpha(vu, u^{-1})\!=\!%
(u^{-1}vu, u^{-1})\!=\!(v^{u}, u^{-1})
$$
$s'(u, v)\!=\!\beta\alpha(vu, v)\!=\!\beta(u^{-1}\!vu, u^{-1}\!v)\!=\!%
(v^{vu}, u^{-1}\!v^{-1}\!v)\!=\!(v^{u}, u^{-1}\!)\!=\!s(u, v)$%
, ainsi~:
$$s^{2}(u, v) = s(v^{u}, u^{-1}) = ((u^{-1})^{v^{u}}, (v^{-1})^{u})
= ((u^{-1})^{vu}, (v^{-1})^{u})
= (u^{-1}, v^{-1})^{vu}
$$
$$s^{4}(u, v)\!=\!s^{2}((u^{-1}, v^{-1})^{vu})\!=\!(u, v)^{vu(v^{-1}u^{-1})^{vu}}\!=\!%
(u, v)^{[v,\, u]}\!=\!\varphi_{[u, v]}(u, v)
$$
\hfill\carre

Comme%
\footnote{car
$\rho(\alpha) = \begin{pmatrix}1& -1\cr 0& 1\cr\end{pmatrix},\,
\rho(\beta) =\begin{pmatrix}1& 0\cr 1& 1\cr\end{pmatrix}\in SL(2, {\Bbb Z})$%
, donc
$\Psi(B_{3})\subset Aut^{+}(F_{2})$%
.
}\
$\ker(\Psi)\!\subset\!ker(\rho\circ\Psi)%
\!\subset\!ker(\rho^{+}\circ\Psi)\!=\,\langle s^{4}\rangle$%
, celui de la pr\'esen\-ta\-tion d'Artin
$({\goth A})$%
, l'int\'eriomorphisme 
$s^{4} = \varphi_{[u, v]}$
\'etant d'ordre infini,
$\Psi$
est injectif.\break\null\hfill\carre\findem

D\'esormais on renote la paire de g\'en\'erateurs de
$F_{2}$ 
en
$t = (t_{1}, t_{2}) := (u, v)$%
.

\Defnns Soit un mot non vide
$W = W(T) = T_{i_{0}}^{\epsilon_{0}}\!\cdots T_{i_{k}}^{\epsilon_{k}}\cdots T_{i_{n}}^{\epsilon_{n}}$%
,
$(i_{j}, \epsilon_{j})\!\in\!\{1, 2\}\!\times\!\{-1, 1\}$
 en les ind\'etermin\'ees%
\footnote{non-commutatives, on notera
$T = (T_{1}, T_{2})$
la paire de ces ind\'etermin\'ees.
}\
$T_{1}, T_{2}$
et leurs {\it inverses\/}
$T_{1}^{-1}, T_{2}^{-1}$%
, son {\it mot inverse\/} est
$W^{-1}\!=\!W^{-1}(T)\!=\:T_{i_{n}}^{-\epsilon_{n}}\cdots T_{i_{k}}^{-\epsilon_{k}}\cdots T_{i_{0}}^{-\epsilon_{0}}$%
. Le mot vide est not\'e
$1:=\emptyset$%
.

Si
$k>0, i_{k-1}\!=\!i_{k},\ \epsilon_{k-1}+\epsilon_{k}\!=\!0$
et
$\{j_{0},\ldots, j_{n-2}\}\!\!=\!\!%
\{0,\ldots, n\}\!\setminus\!\{k-1, k\}$
, dans le m\^eme ordre,
$W'\!=\!W_{k, k+1}\!=\!T_{i_{j_{0}}}^{\epsilon_{j_{0}}}\!\cdots T_{i_{j_{n-2}}}^{\epsilon_{j_{n-2}}}$
est dit {\it r\'eduction \'el\'ementaire\/} de
$W$%
.

On notera plus visuellement si%
\footnote{et
$W_{0, 1}  = \{T_{i_{0}}^{\epsilon_{0}}T_{i_{1}}^{\epsilon_{1}}\}T_{i_{2}}\!\cdots T_{i_{n}}^{\epsilon_{n}},
W_{n-1, n} = T_{i_{0}}^{\epsilon_{0}}\!\cdots T_{i_{n-2}}^{\epsilon_{n-2}}\{T_{i_{n-1}}T_{i_{n}}^{\epsilon_{n}}\}$%
.
}\
$1 <k <n,\ W_{k-1, k} = T_{i_{0}}^{\epsilon_{0}}\cdots 
\{T_{i_{k-1}}^{\epsilon_{k-1}}T_{i_{k}}^{\epsilon_{k}}\}\cdots T_{i_{n}}^{\epsilon_{n}}$%
.

Une {\it r\'eduction\/} de
$W$
est un mot
$W''$
tel qu'il y a
$W=W_{0},\cdots,W_{p} = W''$%
une suite de mots telle que pour
$1\leq k\leq p, W_{k}$
est r\'eduction \'el\'ementaire de
$W_{k-1}$%
.

Si un mot
$W=ABC$
est juxtapos\'e de trois mots
$A, B, C$
avec
$A$
(resp.
$B, C$%
)
se r\'eduisant sur
$1$%
, on notera
$W=\{A\}BC$
(resp.
$A\{B\}C, AB\{C\}$%
).

Un mot
$W,$
est {\it \'equivalent\/} \`a 
$W'$%
, not\'e
$W\approx W'$
s'il y a
$W\!=\!W_{0},\ldots, W_{r}\!=\!W'$
avec
 $W_{l}$
 r\'educ\-tion \'el\'ementaire de
 $W_{l-1}$
 ou 
 $W_{l+1}$
 (supposant
 $l>0,\ l<r$
 dans le premier et le second cas).
 
 Un mot
$W = T_{i_{0}}^{\epsilon_{0}}\!\cdots T_{i_{n}}^{\epsilon_{n}}$
est dit {\it r\'eduit\/} s'il n'a pas de r\'eduction \'el\'ementaire.

Un mot
$W = T_{i_{0}}^{\epsilon_{0}}\!\cdots T_{i_{n}}^{\epsilon_{n}}$
est {\it cycliquement r\'eduit\/} s'il est r\'eduit et si
$T_{i_{0}}^{\epsilon_{0}}\ne T_{i_{n}}^{-\epsilon_{n}}$%
.

Pour tout mot
$W$
il y a un unique%
\footnote{sans unicit\'e de la r\'eduction
$\{II^{-1}\}I$
et
$I\{I^{-1}I\}$
sont deux r\'eductions de
$II^{-1}I$
\`a
$\widetilde{II^{-1}I} = I$%
. Pour une d\'emonstration, voir Th\'eor\`eme 1 et Lemme 1 de \cite{Se}.
}\
mot r\'eduit, not\'e
$\widetilde{W}$%
, r\'eduction de
$W$%
. Cette {\it r\'eduction r\'eduite\/} de
$W$
ne d\'epend que de la classe d'\'equivalence de
$W$%
.

Le nombre
$n+1$
de termes de
$W = T_{i_{0}}^{\epsilon_{0}}\cdots T_{i_{n}}^{\epsilon_{n}}$
 est la
{\it longueur\/}
$l_{T}(W) = l(W)$
du mot
$W$%
.

En substituant les g\'en\'erateurs
$t_{1}, t_{2}$
de 
$F_{2}$
dans les ind\'etermin\'ees
$T_{1}, T_{2}$%
, on a les m\^emes d\'efinitions et notations pour les \'el\'ements du
groupe libre \`a deux g\'en\'erateurs~:
 $$w = t_{i_{0}}^{\epsilon_{0}}\!\cdots t_{i_{n}}^{\epsilon_{n}}\!\in\!F_{2},\
 (i_{j}, \epsilon_{j})\!\in\!\{1, 2\}\!\times\!\{-1, 1\}$$
\finc
\Defn
Un endomorphisme 
$x_{0} : F_{2}\!\rightarrow\!F_{2}$
induit un endomorphisme~:
$$\overline{x_{0}}\!=\!ab(x) : {\Bbb Z}^{2}\!=\!%
ab(F_{2}) \rightarrow ab(F_{2})\!=\!{\Bbb Z}^{2}$$
%
 de l'ab\'elianis\'e
 ${\Bbb Z}^{2}$%
 . On dira que l'endomorphisme
$x_{0}$
est {\it direct\/} si 
$det(\overline{x_{0}})>0$%
.

\finc
\Thc Th\'eor\`eme  {\rm (de Nielsen, 1917)}|
${\romannumeral 1})$
Un endomorphisme surjectif
$x$
du groupe libre \`a deux g\'en\'erateurs
$F_{2}$
est un automorphisme.

${\romannumeral 2)}$
Le groupe d'automorphismes directs du groupe libre
$F_{2}$
est produit semi-direct suivant
$\psi$
amalgam\'e sur le sous-groupe
$Y$
d'indice
$2$
dans le centre de
$B_{3}$%
~:
$$Aut^{+}(F_{2}) = Int(F_{2})\rtimes_{<\varphi_{[t_{1}, t_{2}]}\,=\,s^{4}>}B_{3}\ .$$

${\romannumeral 3)}$
Le groupe
$Int(F_{2})$
des automorphismes int\'erieurs de
$F_{2}$ 
est le noyau de~:
$$\rho: Aut(F_{2})\rightarrow Aut({\Bbb Z}^{2})=GL(2, {\Bbb Z})\ .$$

${\romannumeral 4)}$
Un endomorphisme 
 $\varphi:\!F_{2}\!\rightarrow\!F_{2},\ \varphi(t_{i})\!=\!X_{i},\ i\!=\!1, 2$
du groupe libre \`a deux g\'en\'erateurs
est un automorphisme direct de ce groupe libre
si et seulement si il pr\'eserve \`a conjugaison pr\`es
le commutateur
$[t_{1}, t_{2}]$
des g\'en\'erateurs~:$$[X_{1}, X_{2}]=X_{1}^{-1}X_{2}^{-1}X_{1}X_{2}=[t_{1}, t_{2}]^{w}=w^{-1}[t_{1}, t_{2}]w
\ .\leqno{(Cw)}$$
De plus, il est dans
$B_{3}$
si et seulement si
$[X_{1}, X_{2}]\!=\![t_{1} ,t_{2}]$%
, c.\`a.d. on a 
$(C_{1})$%
.
\finc

\Rmc {\petcap Remarques\/}|
$(1)$
${\romannumeral 1)}$
vaut pour tout groupe libre de type fini%
\footnote{nous donnerons  cependant ici ausi une d\'emonstration de
${\romannumeral 1)\/}$
\og\`a la Nielsen\fgf.
}\
$F$
: 
le gradu\'e central descendant de
$F$
\'etant engendr\'e par son terme de degr\'e
$1$%
, un endomorphisme surjectif
$x$
de
$F$
est surjectif au niveau du gradu\'e central descendant
$gr_{\goth c}(F)$
donc, puisque chaque terme est de dimension finie, injectif et,
puisqu'un \'el\'ement non nul de
$F$
se d\'etecte dans
$gr_{\goth c}(F),\ x$
est aussi injectif donc un isomorphisme.
\hfill\carre

$(2)$
Le point
${\romannumeral 3)}$
est particulier \`a la dimension
$2$%
~:
$$\psi : F_{3}=\langle t_{1}, t_{2}, t_{3}\rangle\rightarrow F_{3},\
\psi(t_{1}) = t_{1}^{t_{2}}=t_{2}^{-1}t_{1}t_{2},\
\psi(t_{2})=t_{2},\, \psi(t_{3})=t_{3}$$
est dans le noyau de
$\rho_{3}: Aut(F_{3})\rightarrow GL(3, {\Bbb Z})$
mais n'est pas int\'erieur%
\footnote{si
$\psi=\varphi_{w}$
alors
$w$
doit commuter \`a
$t_{2}$%
, donc
$w=t_{2}^{n}$
pour un
$n\in{\bf Z}$%
, mais il doit aussi commuter \`a
$t_{3}$%
, ainsi
$n=0$
et
$\varphi_{w}=\varphi_{1}=Id\ne\psi$
car
$\psi(t_{1})\ne t_{1}$%
.\hfill\carre
}.\hfill\carre

$(3)$
Les \'el\'ements de
$H$
et
$B_{3}$
v\'erifient les conditions%
\footnote{Comme
$H\!:=\!Int(F_{2})\!\cap\!B_{3}$
et
$t\mapsto t^{w}$
est un morphisme, l'isomorphisme int\'erieur
$\varphi_{w^{-1}}$%
, on a
$[t_{1}^{w}, t_{2}^{w}]\!=\![t_{1}, t_{2}]^{w}$
et 
$[\alpha(t_{1}), \alpha(t_{2})]=[t_{1}, t_{1}^{-1}t_{2}]=
t_{1}^{-1}(t_{1}^{-1}t_{2})^{-1}t_{1}t_{1}^{-1}t_{2}=
t_{1}^{-1}(t_{2}^{-1}t_{1})t_{1}t_{1}^{-1}t_{2}=
[t_{1}, t_{2}]$
et~: 
$[\beta(t_{1}), \beta(t_{2})]=[t_{2}t_{1}, t_{2}]=
(t_{2}t_{1})^{-1}t_{2}^{-1}(t_{2}t_{1})t_{2}=
t_{1}^{-1}t_{2}^{-1}t_{2}^{-1}t_{2}t_{1}t_{2}=
[t_{1}, t_{2}]$%
.\hfill\carre
}
$(C_{w})$
et
$(C_{1})$%
.\hfill\carre

$(4)$
Si
$x$
v\'erifie
$(C_{1})$%
, alors
$x\!=\!\varphi_{w}b$
o\`u
$w\!\in\!F_{2},\ b\!\in\!B_{3}$
et, si
$^{w}[t_{1}, t_{2}]\!=\![t_{1}, t_{2}]$%
, alors
$w\!=\![t_{1}, t_{2}]^{n}\!\in\!\langle [t_{1}, t_{2}]\rangle$%
, ainsi le \og De plus%
$\cdots$%
\fg 
suit%
\footnote{En effet, comme
$\varphi_{[t_{1}, t_{2}]^{n}}=s^{4n}$%
, on a
$x=s^{4n}b\in B_{3}$%
.\hfill\carre}%
$\,$
de
${\romannumeral 4})$
et
${\romannumeral 2})$%
.\hfill\carre
\finc
Dans la suite on dira
\og Nielsen
$\ast$%
\fg
pour
\og le point
$\ast$
du {\petcap Th\'eor\`eme\/} de Nielsen\fgf.

\vskip2mm

\Defns
La {\it taille\/}
d'un endomorphisme
$x$
d\'efini
par la paire de mots r\'eduits%
\footnote{en la paire d'ind\'etermin\'ees
$T=(T_{1}, T_{2})$%
.
}~:\hfill\break
$X\!=\!(X_{1}, X_{2}),\, x(t_{i})\!=\!X_{i}(t),\, i\!=\!1, 2$
est
$|x|:=l_{T}(X_{1})+ l_{T}(X_{2})$%
. L'endomorphisme
$x$
est
{\it minimal\/}
si sa taille
$|x|$
est minimale dans la double classe
$HxH$%
.

Soit
$M\!=\!X_{i_{0}}^{\epsilon_{0}}\!\!\cdots\!X_{i_{n}}^{\epsilon_{n}}$
un mot r\'eduit en une paire d'ind\'etermin\'ees not\'ee abusivement
$X\!=\!(X_{1}, X_{2})$
et
$N\!=\!N(T)$%
, le mot (non r\'eduit), obtenu apr\`es avoir substitu\'e dans
$X$
la paire 
$X\!=\!(X_{1}(T), X_{2}(T))$
d\'efinissant 
$x$%
.

Il y a des mots r\'eduits%
\footnote{en la paire d'ind\'etermin\'ees
$T$%
, certains de ces mots
$C_{i}, D_{j}$
pouvant \^etre vides.
}\
$C_{0}, D_{0}, \cdots, C_{k}, D_{k},\cdots, C_{n-1}, D_{n-1}, C_{n}$
tels que~:
$$X_{i_{0}}^{\epsilon_{0}}\!=\!C_{0}D_{0},\ \hbox{\rm pour\/}\
0\!<\!k\!<\!n,\ X_{i_{k}}^{\epsilon_{k}}\!=\!D_{k-1}^{-1}C_{k}D_{k},\ \hbox{\rm et \/}\
X_{i_{n}}^{\epsilon_{n}}\!=\!D_{n-1}^{-1}C_{n}$$
 comme mots r\'eduits et, pour
$0<k\leq n$%
, le mot
$C_{k-1}C_{k}$
est r\'eduit, ils peuvent \^ etre d\'efinis inductivement dans cet ordre
 en partant de la gauche%
\footnote{$C_{0}\{D_{0}\!D_{0}^{-1}\!\}E_{1}\!X_{i_{2}}^{\epsilon_{2}}\!\cdots\!%
X_{i_{n}}^{\epsilon_{n}}$
est la plus grande suite de r\'eductions \'el\'ementaires entre
$X_{i_{0}}^{\epsilon_{0}}$
et
$X_{i_{1}}^{\epsilon_{1}}$
dans
$X_{i_{0}}^{\epsilon_{0}}\cdots X_{i_{n}}^{\epsilon_{n}}$%
, puis
$C_{0}C_{1}\{D_{1}\!D_{1}^{-1}\!\}E_{2}\!X_{i_{3}}^{\epsilon_{3}}\!\cdots\!%
X_{i_{n}}^{\epsilon_{n}}$
est la plus grande suite de r\'eductions \'el\'ementaires entre
$E_{1}$
et
$X_{i_{2}}^{\epsilon_{2}}$
dans
$C_{0}E_{1}X_{i_{2}}^{\epsilon_{2}}\cdots
X_{i_{n}}^{\epsilon_{n}}$%
\dots
}\
ou de la droite.

Un endomorphisme direct est dit
{\it admissible\/}
s'il est surjectif ou satisfait
$(C_{w})$%
.

\finc

\Thc {\sc Lemme.}|
Si l'endomorphisme
$x$
est admissible minimal alors pour
$0\!\leq\!i\!\leq\!n$
le mot
$C_{i}\ne\emptyset$ 
est non vide  donc
$C_{0}\cdots C_{n}$
est la  r\'eduction r\'eduite
$\widetilde{N}$
du mot
$N$%
.
\finc

Les
$C_{i},\, i = 0,\ldots, n$
(non vides, d'apr\`es ce {\petcap Lemme\/}) sont dits {\it subsistants\/} de
$N$%
.
\vskip2mm
\Demd{u Th\'eor\`eme}|%
${\romannumeral 3)}$
Comme 
$Int(F_{2})\!\subset\!ker(\rho)\!\!\subset\!\!Aut^{+}(F_{2})$
et, d'apr\`es la pr\'esentation
$({\goth A\/})$%
, on a
$B_{3}\!\cap\!ker(\rho)\!=\!Y$%
, le point
${\romannumeral 3)}$
suit de
${\romannumeral 2)}$%
.\hfill\carre

${\romannumeral 1)}$
et
${\romannumeral 2)}$
Soit
$x_{0}$
admissible et
$x\!=\!k\!\circ\!x_{0}\!\circ\!h,\, h,k\!\in\!H$%
, minimal pour la taille
$|x|$
de
$x$
dans la double classe
$Hx_{0}H$
de
$x_{0}$
suivant
$H$%
.

Si
$\tau\!\in\!Aut(F_{2}), \tau(t_{i})\!=\!t_{3-i},\, i\!=\!1, 2$
alors
$det(\overline{\tau})\!=\!-1$
et
$|x\tau|\!=\!|x|$
donc, quitte \`a prendre
$x\tau$%
, il suffit  dans
${\romannumeral 1)}$%
,  de consid\'erer un endomorphisme surjectif direct.

Par surjectivit\'e de
$x$
il y a une paire%
\footnote{de mots r\'eduits en la paire d'ind\'eter\-min\'ees
$X = (X_{1}, X_{2})$%
.
}\
$Y\!=\!(Y_{1}, Y_{2})$
telle que
$\widetilde{Y_{k}}(T)\!=\!T_{k}$
pour
$k=1, 2$%
. Ainsi, par le {\petcap Lemme\/},
$Y_{1}$
et
$Y_{2}$
sont de longueur
$1$%
, donc aussi
$X_{1}$
et
$X_{2}$%
.

L'endomorphisme admissible minimal
$x$
\'etant direct on a~:
$$\hbox{\rm soit\/}\ X = T,\ x(t) = (t_{1}, t_{2}) = t$$
$$\hbox{\rm soit\/}\ X = (T_{1}^{-1}, T_{2}^{-1}),\ x(t) =
(t_{1}^{-1}, t_{2}^{-1}) = s^{2}(t)^{t_{1}^{-1}t_{2}^{-1}}: x = \varphi_{t_{2}t_{1}}s^{2}\in H$$
$$\hbox{\rm soit\/}\ X = (T_{2}^{\epsilon}, T_{1}^{-\epsilon}),\ x(t) =
(t_{2}^{\epsilon}, t_{1}^{-\epsilon}) =
\bigl(s^{\epsilon}(t)^{t_{1}^{-1}}\bigr): x = \varphi_{t_{1}}s^{\epsilon}\in H,
\ \epsilon = -1, 1.$$

Dans tous ces cas l'endomorphisme admissible minimal
$x\!=\!kx_{0}h$
 est dans le sous-groupe
$H$%
, donc
$x_{0}$%
, est dans 
$H\!=\!Int(F_{2})\rtimes_{<\varphi_{[t_{1}, t_{2}]}=s^{4}>}B_{3}$%
, d'o\`u
${\romannumeral 1)}$
et
${\romannumeral 2)}$%
.
\hfill\carre

${\romannumeral 4)\/}$
D'apr\`es
la {\petcap Remarque\/}
$(3)$
on peut demander \`a l'endomorphisme
$x$%
, parmi ceux minimisant la taille
$|x|$
dans la double classe
$Hx_{0}H$%
, de minimiser aussi
$l(w)$%
. Alors%
\footnote{Sinon pour un
$i\!\in\!\{1, 2\}, w\!=\!t_{i}^{-1}w', l(w')\!=\!l(w)-1$
et, comme pour
$ i\!=\!1, 2$
on a~:
$
[t_{1}, t_{2}]^{t_{i}^{-1}}\!\!=\![t_{2}^{3-2i}, t_{1}^{2i-3}]\!=\!%
[k^{-1}(t_{1}), k^{-1}(t_{2})]
$%
, o\`u
$k\!=\!k_{i}\!=\!(s^{2i-3})^{t_{1}^{-1}}$
et
$x'\!=\!kx, w''\!=\!k(w')$%
, on a~:
$x'\!\in\!Hx\!\subset\!HxH,\, |x'|\!=\!|x|,\, l(w'')\!=\!l(w')\!<\!l(w)$
et
$[x'(t_{1}), x'(t_{2})]\!=\!k([x(t_{1}), x(t_{2})]\!=\!k([t_{1}, t_{2}]^{t_{i}^{-1}w'})%
\!=\!k([k^{-1}(t_{1}), k^{-1}(t_{2})]^{w'})\!=\!([kk^{-1}(t_{1}), kk^{-1}(t_{2})]^{k(w')}%
\!=\![t_{1}, t_{2}]^{w''}$%
, contredi\-sant que
$x$
minimise
$|x|$%
, puis
$l(w)$
dans la double classe
$HxH$%
.\hfill\carre
}\
$l(t_{2}w) = l(w)+1 = l(t_{1}w)$
et
 $w^{-1}[t_{1}, t_{2}]w$
est r\'eduction r\'eduite de
$[X_{1}, X_{2}]$%
.

Si
$w\!\ne\!1$%
, donc
$w\!=\!w'a,\, a\!=\!t_{k}^{\pm 1}$%
, alors, par le {\petcap Lemme\/},
$X_{i}\!=\!X''_{i}a$
pour
$i\!=\!1,2$
et
$x'\!=\!\varphi_{a}x=x^{a^{-1}}$%
, d\'efini par
$(X'_{1}, X'_{2})=(\widetilde{aX''}_{1}, \widetilde{aX''_{2}})$%
, v\'erifie
$[X'_{1}, X'_{2}]\!=\![t_{1}, t_{2}]^{w'}$
et
$|x'|\leq |x|$
contredisant la minimalit\'e de
$|x|$%
, puis
$l(w)$
dans la double classe
$Hx_{0}H$%
.

On a donc
$(C1)\ [X_{1}, X_{2}]=[t_{1}, t_{2}]$
et le {\petcap Lemme\/} donne les \'egalit\'es r\'eduites~:
$$X_{1}^{-1}=t_{1}^{-1}D_{0},\ X_{2}^{-1}=D_{0}^{-1}t_{2}^{-1}D_{1},\
X_{1}=D_{1}^{-1}t_{1}D_{2},
\ X_{2}=D_{2}^{-1}t_{2}$$

On a
$l(X_{i}^{-1})\!=\!l(X_{i}),\, i\!=\!1,2$
donc
$l(D_{0})\!=\!l(D_{1})+l(D_{2})\!=\!l(D_{1})+l(D_{0})+l(D_{1})$%
, ainsi
$l(D_{1})\!=\!0$
et
$X_{1}^{-1}\!=\!t_{1}^{-1}D_{0}, X_{2}^{-1}\!=\!D_{0}^{-1}t_{2}^{-1},\, X_{2}X_{1}\!=\!(t_{2}D_{0})(D_{0}^{-1}t_{1})\!=\!t_{2}t_{1}$%
.

Comme
$l(X_{1})+l(X_{2})\!=\!|x|\!\leq\!|x\beta|\!=\!l(\widetilde{X_{2}X_{1}})+l(X_{2})\!=\!2+l(X_{2})$%
, on a
$l(D_{0})\!\leq\!1$%
.

Ainsi
$l(D_{2})\!\leq\!1$%
. Les
$X_{i}$
\'etant r\'eduits
$t_{1}^{-1}\!\ne\! D_{2}\!\ne\! t_{2}$%
. Or%
\footnote{sinon
$X_{1}\!=\!t_{1}^{2}$
ou
$X_{2}\!=\!t_{2}^{2}$
et, dans le terme de degr\'e
$2$
du gradu\'e central descendant
$[F_{2}, F_{2}]/\bigl[[F_{2}, F_{2}], F_{2}\bigr]\simeq{\Bbb Z}$%
, la classe du commutateur
$[X_{1}, X_{2}]$
serait divisible par
$2$
contredisant que, d'apr\`es la condition
$(C1)$%
, elle est celle de son g\'en\'erateur.
}\
$D_{2}\!\in\!\{t_{1}, t_{2}^{-1}\}$%
, donc
$D_{2}\!=\!1$
et
$X_{1}\!=\!t_{1}, X_{2}\!=\!t_{2},\, x\!=\!1$%
, ainsi
$x_{0}\!=\!k^{-1}xh^{-1}\!=\!k^{-1}h^{-1}\!\in\!H$%
.\hfill\carre\findem

\Demd{u Lemme}|
Rappelons que
$(X_{1}, X_{2})$
est la paire de mots en la paire
$T\!=\!(T_{1}, T_{2})$
d'ind\'etermin\'ees telle que l'endomorphisme
$x$
est donn\'e par
$x(t_{i})\!=\!X_{i}(t), i\!=\! 1, 2$%
.
Si, pour un
$1\leq k\leq n,\, i_{k}\ne i_{k-1}$%
, on note~:
$$U=X_{i_{k-1}}^{\epsilon_{k-1}},\ V=X_{i_{k}}^{\epsilon_{k}},\ \nu=\epsilon_{k-1},\
\epsilon=\epsilon_{k-1}\epsilon_{k}\ .$$

\finc

\Thc 
{\sc Affirmation 1}|
Il y a
$h\in H$
tel que
$xh$
minimise encore la taille dans la double classe
$Hx_{0}H$
et est d\'efini par l'un des deux cas~:

{\romannumeral 1}) (cas I)\quad
$(U, V)$

{\romannumeral 2}) (cas II)\quad
$(V, U)$%
\ .
\finc

\Dem
Quitte%
\footnote{puisque
$x\varphi_{t_{1}}s\!\in\!xH\!\subset\!Hx_{0}H$%
, \'etant d\'efini par
$t\mapsto (X_{2}, X_{1}^{-1})$%
, a m\^eme taille que
$x$%
.
}\
\`a changer
$x$
en
$x\varphi_{t_{1}}s$%
, on peut supposer
$\epsilon = +1$%
.\hfill\break
Alors, si
$\nu = -1,\ h = \varphi_{t_{2}t_{1}}s^{2}$
convient et si
$\nu = 1,\ \{U, V\} = \{X_{1}, X_{2}\}$%
. \hfill\findem

Par abus, on renomme d\'esormais le
$xh$
donn\'e par l'{\petcap affirmation 1\/} en
$x$%
.

\finc

\Thc {\sc Compl\'ement 1}|
Suivant le cas I ou II, les paires de mots r\'eduits~:
$$(U, \widetilde{UV}),\ (U, \widetilde{VU}),\
(\widetilde{VU}, V)\leqno{(\hbox{\rm cas I\/})}$$
$$(\widetilde{UV}, U),\ (\widetilde{VU}, U),\
(V, \widetilde{VU})\leqno{(\hbox{\rm cas II})}$$
d\'efinissent des endomorphismes  de la double classe
$Hx_{0}H$
de
$x_{0}$%
, \`a savoir~:
$$x\alpha^{-1},\, \varphi_{U}^{-1}x\alpha^{-1},\, x\beta\quad
\hbox{\rm resp.\/}\quad x\beta,\, \varphi_{U}^{-1}x\beta,\,
x\alpha^{-1} \quad\quad\quad\quad\carre$$

\finc
\Demo{pour les subsistants extr\^emes
$C_{0}, C_{n}$%
}|
Si
$C_{0} = \emptyset$
est vide on a
$X_{i_{0}}^{\epsilon_{0}}=D_{0}$
et, comme le mot
$M$
est r\'eduit,
$i_{1}\ne i_{0}$%
. 

Alors%
\footnote{avec les notations ci-dessus avec
$k\!=\!1$%
.
}\
$l(\widetilde{UV})=l(V)-l(U)<l(V)$
et suivant le cas I ou II, d'apr\`es le compl\'ement de l'{\petcap Affirmation 1\/},
l'endomorphisme
$x\alpha^{-1}$
ou
$x\beta$
contredit que
$xh$
minimise la taille dans la double classe
$Hx_{0}H$%
.\hfill\carre

Ainsi le subsistant initial
$C_{0}\ne\emptyset$
est non vide. En consid\'erant le mot inverse
$M^{-1}\!=\!X_{i_{n}}^{-\epsilon_{n}}\cdots X_{i_{0}}^{-\epsilon_{0}}$
on \'etablit de m\^ eme que
$C_{n}\ne\emptyset$
est non vide.\hfill\carre
\vskip2mm
\Demo{pour les subsistant int\'erieurs}|
Supposons, pour un
$0\!<\!k<n,$
$C_{k}\!=\!\emptyset$
 vide, alors
$X_{i_{k}}^{\epsilon_{k}}\!=\!D_{k-1}^{-1}D_{k}$
et
$D_{k}, D_{k-1}^{-1}$%
, ne pouvant \^etre subsitant final, ou inital, vide du mot 
$M_{k-}\!=\!X_{i_{0}}\cdots X_{i_{k}}$%
ou 
$M_{k+}\!=\!X_{i_{k}}\cdots X_{i_{n}}$%
, on a~:
$D_{k}\!\ne\!\emptyset\!\ne\!D_{k-1}$%
.

\finc

\Thc {\sc Affirmation 2}|
$i_{k-1}\ne i_{k}\ne i_{k+1}$
(donc%
\footnote{c'est ici que 
$x\!\in\!End(F_{2})$%
, plut\^ot que
$End(F)$
(%
$F$
un groupe libre de type fini), est crucial!
}\
$i_{k-1}=i_{k+1}$%
)
et
$\epsilon_{k-1} = \epsilon_{k+1}$%
.
\finc
\vskip2mm
\Dem
Si
$i_{k-1}\!=\!i_{k}\!\!=\!i_{k+1}$%
, alors, puisque
$M$
est r\'eduit, on a~:\hfill\break
$\epsilon_{k-1}\!=\!\epsilon_{k}\!=\!\epsilon_{k+1}$
et\nobreak%
 \footnote{avec les conventions
$D_{-1} = \emptyset = D_{n + 1}$%
.
}\
$D_{k-2}^{-1}C_{k-1}D_{k-1} = D_{k-1}^{-1}D_{k} =  D_{k}^{-1}C_{k+1}D_{k+1}$%
.

Suivant que
$l(D_{k-1})\!\leq\!l(D_{k})$
ou
$l(D_{k-1})\!\geq\!l(D_{k})$%
, on a
$X_{i_{k}}^{\epsilon_{k}}\!=\!D_{k-1}^{-1}ED_{k-1}$
ou
$X_{i_{k}}^{\epsilon_{k}}\!=\!D_{k}^{-1}ED_{k}$%
, on notera selon le cas
$D=D_{k-1}$
ou
$D_{k}$%
.

Comme%
\footnote{par l'hypoth\`ese faite sur
$x_{0}$
 (celle de
${\romannumeral 1)\/}$
ou
${\romannumeral 4)\/}$
du th\'eor\`eme).
}\
on a
$X_{i_k}^{\epsilon_{k}}\ne 1$
donc
$E\ne\emptyset,\ X_{i_{k-1}}^{\epsilon_{k-1}}X_{i_{k}}^{\epsilon_{k}}X_{i_{k+1}}^{\epsilon_{k+1}}=D^{-1}E^{3}D$
et le terme
$E$ central, celui qui est dans
$X_{i_{k}}^{\epsilon_{k}}=D^{-1}ED$%
, ne peut dispara\^{\i}tre.\hfill\carre

Si
$i_{k-1}\!=\!i_{k}\!\ne\!i_{k+1}$
ou
$i_{k-1}\!\ne\!i_{k}\!=\!i_{k+1}$%
, quitte \`a consid\'erer le mot inverse
$M^{-1}$%
, on supposera que l'on est dans le second cas
$i_{k}\!=\!i_{k+1}$
alors, avec les notations pr\'ec\'edant l'{\petcap Affirmation 1\/}
et car
$M$
est r\'eduit, on a
$X_{i_{k-1}}^{\epsilon_{k-1}}X_{i_{k}}^{\epsilon_{k}}X_{i_{k+1}}^{\epsilon_{k+1}}\!=\!UVV$%
.

Comme
$\emptyset\!\ne\!D_{k},\ D_{k}\!=\!D'_{k}a,\, a\!=\!t_{i}^{\pm 1}$
derni\`ere lettre du
$V$
central donc
$a^{-1}$
est  premi\`ere du
$V$
de droite et%
\footnote{puisque
$D_{k-1}^{-1}$
est segment initial du
$V$
central.
}\
de
$D_{k-1}^{-1}$
donc
$a$
est aussi derni\`ere lettre de
$U$%
.

Ainsi
$l(V^{a^{-1}})\!=\!l(V)-2\!<\!l(V),\ l(U^{a^{-1}})\!\leq\!l(U)$
et
$\varphi_{a}x\!\in\!Hx_{0}H$
contredit que 
$|x|$
minimise la taille dans la double classe
$Hx_{0}H$%
.\hfill\carre

Si
$\epsilon_{k-1}\!=\!-\epsilon_{k+1}$
on a
$X_{i_{k-1}}^{\epsilon_{k-1}}X_{i_{k}}^{\epsilon_{k}}X_{i_{k+1}1}^{\epsilon_{k+1}}\!=\!UVU^{-1}$%
, donc
$D_{k-1}$
et
$D_{k}$
sont segments finaux de
$U\!=\!(U^{-1})^{-1}$%
, contredisant%
\footnote{puisque
$D_{k-1}\!\ne\!\emptyset\!\ne\!D_{k}$%
.
}\
que
$V\!=\!D_{k-1}^{-1}D_{k}$
est r\'eduit.\hfill\findem

Ainsi, posant
$D = D_{k-1}^{-1},\, E = D_{k}$%
, on a%
\footnote{avec les notations
$U\!=\!X_{i_{k-1}}^{\epsilon_{k-1}}, V\!=\!X_{i_{k}}^{\epsilon_{k}}$
de l'Affirmation 1 car,
par l'Affirmation 2, 
$i_{k-1}\ne i_{k}$%
.
}\
$V = DE$
et
$U_{-}D^{-1} = U = E^{-1}U_{+}$%
.

Quitte \`a prendre l'inverse
$(UVU)^{-1} =
U^{-1}V^{-1}U^{-1}$
du mot
$UVU$
et \`a renommer
$(U, V, D, E):=(U^{-1}, V^{-1}, E^{-1}, D^{-1})$%
, on supposera
$l(D)\leq l(E)$%
.

\parindent=0pt\par{\sl Fin de la d\'emonstration du {\petcap Lemme\/}}:%
\parindent=20pt

Si
$l(U_{+})\!<\!l(E)$
(donc
$l(\widetilde{VU})\!<\!l(V)$%
)%
, suivant le cas I ou II,
$\varphi_{U}^{-1}x\alpha^{-1}\!$
ou
$\varphi_{U}^{-1}x\beta\!$
contredit encore que
$x$
minimise la taille dans la double classe
$Hx_{0}H$%
.
\hfill\carre

Ainsi
$l(U_{+})\geq l(E)\geq l(D)$
et il y a un mot, \'eventuellement vide,
$F$
tel que
$U_{+}$
a  l'expression r\'eduite
$U_{+} = FD^{-1}$
et donc
$U= E^{-1}FD^{-1}$%
, ainsi~:
$$\varphi_{D}^{-1}(V)=(DE)^{D} = ED\
\hbox{\rm et\/}\
\varphi_{D}^{-1}(U) = (E^{-1}FD^{-1})^{D} = D^{-1}E^{-1}F$$
donc
$\varphi_{D}^{-1}(\widetilde{VU}) = \varphi_{D}^{-1}(V)\varphi_{D}^{-1}(U) = F$
et
$l(\varphi_{D}^{-1}(V))=l(V),\ l(\varphi_{D}^{-1}(\widetilde{VU}))<l(U)$%
.

Donc,  suivant le cas I ou II, que
$x$
minimise la taille dans
$Hx_{0}H$
est encore contredit par l'endomorphisme surjectif
$\varphi_{D}^{-1}x\beta$
ou
$\varphi_{D}^{-1}x\alpha^{-1}$%
.

Aucun subsistant
$C_{k}$ ne pouvant ainsi \^etre vide, le Lemme est \'etabli.\hfill\carre\findem


\section*{APPENDICE C\\
Forme normale dans
$B_{3}$
et classes de conjugaison de torsion.
}

On utilise les notations du \S {\bf 1\/}~:
$B_{3}=\langle a, b\,;\, aba=bab\rangle$
groupe des tresses \`a trois brins,
d'apr\`es la Proposition de l'appendice B, identifi\'e \`a 
$\Psi(B_{3})\subset Aut^{+}(F_{2})$%
.

Ainsi
$s\!=\!aba\!=\!bab\!: (u, v)\!\mapsto\!(v, u^{-1})^{u},
c\!=\!s^{2}\!=\!ababab, (u, v)\!\mapsto\!(u^{-1}, v^{-1})^{vu}$
et
$d\!=\!c^{2},\, (u, v)\mapsto(u, v)^{[v, u]} $
qui engendrent respectivement le centre 
$Z\!=\!\langle c\rangle$
de
$B_{3}$
et son sous-groupe d'indice
$2,\ Y=\langle d\rangle\,\subset\,\langle c\rangle=Z$%
.\hfill\break
On note
$A
, B, S, C\in SL(2, {\Bbb Z\/})$
leurs images par
$\varrho:=\rho^{+}\!\circ\!\Psi: Aut^{+}(F_{2})\rightarrow SL(2, {\Bbb Z\/})$%
.

\Thc Proposition|
Soit une tresse \`a trois brins
$w\!\in\!B_{3}$%
. Alors il y a~:\hfill\break
$n\!\in\!{\Bbb Z},\, \nu_{1}, \nu_{2}\!\in\!\{0, 1\}$
et
$W(a^{-1}, b)$%
, un mot du mono\"{\i}de engendr\'e par
$a^{-1}$
et
$b$,\
t.q.~:
$$w\!=\!s^{\nu_{1}}W(a^{-1}, b)s^{\nu_{2}}c^{n}$$

De plus, si
$w\!\not\in\!\langle s\rangle$%
, cette \'ecriture est  uniquement%
\footnote{si
$w\!=\!s^{m}\!\in\!Z$
il n'y a pas unicit\'e~:
$c^{n}\!=\!s^{2n}\!=\!s1sc^{n-1}$
(resp.
$s1c^{n}\!=\!s^{2n+1}\!=\!1sc^{n}$%
).
}\
d\'etermin\'ee par
$w$%
.

\finc

\Dem
Il y a
${{1}\over{2}}2^{4}$
choix de signes%
\footnote{comme
$det(\varrho(x))=1>0$
la moiti\'e des choix
$\pm$
pour chacun des quatre coefficients.
}\
pour les coefficients de
$M=\varrho(w)$%
, ainsi il y a
$\nu\!\in\!\{0, 1\}^{3}
$
avec
$M_{\nu}\!=\!S^{-\nu_{1}}\varrho(w)S^{-\nu_{2}}C^{-\nu_{0}}$
\`a coefficients positifs ou nuls.

Par une classique%
\footnote{soit
$M\!=\!\begin{pmatrix}e& f\cr g& h\cr\end{pmatrix}\!, \parallel\!M\!\parallel=\!e\!+\!f\!+\!g\!+\!h$%
. Si
$e\!=\!g$
et
$f\!>\!h$
(resp.
$f\!<\!h$%
)
on note
$M_{1}\!=\!AM$
(resp.
$B^{-1}M$%
). Si
$e\!>\!g$
(resp.
$e\!<\!g$%
)
alors
$f\!\geq\!h$
(resp.
$f\!\leq\!h$%
)
et on pose
$M_{1}\!=\!AM$
(resp.
$B^{-1}M$%
).\hfill\break
Dans tous les cas
$M_{1}$
est \`a coefficients positifs ou nuls avec
$\parallel\!\!M_{1}\!\!\parallel<\parallel\!\!M\!\!\parallel$%
(sinon
$e\!\geq\!g+1, h\!\geq\!f+1$
et
$det(M)\geq (g+1)(f+1)-fg=g+f+1>1$%
, sauf si
$M=Id$%
.\dots)
}\
suite de soustraction de lignes%
\footnote{multiplication \`a gauche par
$A$
ou
$B^{-1}$%
, la note pr\'ec\'edente implique l'unicit\'e du processus.
 }%
, on ram\`ene
$M$
\`a
$\begin{pmatrix}1& 0\cr 0& 1\cr\end{pmatrix}$
et donc
$\varrho(w)$
est le produit 
$W(A^{-1}\!, B)$
en sens inverse des matrices inverses 
$A^{-1}, B$%
.

D'apr\`es 
$({\goth A\/})$
l'\'el\'ement
$w'\!=\!s^{\nu_{1}}W(a^{-1}, b)s^{\nu_{2}}c^{\nu_{0}}$ 
v\'eri\-fiant
$\varrho(w)\!=\!\varrho(w')$
diff\`ere de
$w$
par une puissance paire de
$c~:\ w=w'c^{2m}\!=\!s^{\nu_{1}}W(a^{-1}, b)s^{\nu_{2}}c^{\nu_{0}+2m}$%
.
\hfill\findem
\vskip2mm
\Defn
La
{\it taille normale\/}
de
$w\!\in\!B_{3}$
est l'ordre lexicographique de la  paire
$$t_{n}(w)=(2\nu_{1}+\nu_{2},\, l(W))\ .$$
\Thc Corollaire|
Un \'el\'ement
$w_{0}$
de taille normale minimale dans la classe de conjugaison
$w^{B_{3}}$
de l'\'el\'ement 
$w\in B_{3}$ 
est de taille soit~:

${\romannumeral 1)\/}\quad t_{n}(w_{0})=(0, 0),\ 
\hbox{\rm et\/}\ w_{0}=c^{n}$
est central.

${\romannumeral 2)\/}\quad t_{n}(w_{0})=(0, m),\ m\!\ne\!0\
\hbox{\rm et\/}\ w_{0}=Wc^{n},\, W\!\ne\!1\
\hbox{\rm cycliquement r\'eduit.\/}$

${\romannumeral 3)\/}\quad t_{n}(w_{0})=(1, 0)$
et
$w_{0}=sc^{n}=s^{2n+1}$%
.

${\romannumeral 4)\/}\quad t_{n}(w_{0})=(1, 1),\
\hbox{\rm et soit\/}\ w_{0}=a^{-1}s^{2n+1},\
\hbox{\rm soit\/}\
w_{0}=bs^{2n+1}$%
.\finc
\vskip2mm
\Thc  Compl\'ement|
L'image de
$w_{0}$
dans
$B_{3}/Y=SL(2, {\Bbb Z\/})$
est~:

Dans le cas
${\romannumeral 1)\/},\ \overline{w_{0}}=\overline{c}^{n}$
 d'ordre
$2$
si
$n$
est impair, trivial si
$n$
est pair.

Dans le cas
${\romannumeral 2)\/},\ \overline{w_{0}}=\overline{W}\overline{c}^{n}$
 d'ordre infini.

Dans le cas
${\romannumeral 3)\/},\ \overline{w_{0}}=\overline{s}^{2n+1}$
d'ordre
$4$
et conjugu\'e \`a l'un de
$\overline{s}, \overline{s}^{-1}$%
.

Dans le cas
${\romannumeral 4)\/}$
 d'ordre 
$6$
ou
$3$
et conjugu\'e \`a l'un de
$\overline{b}\overline{a}, (\overline{b}\overline{a})^{-1}\ \hbox{\rm ou\/}\
(\overline{b}\overline{a})^{2},
(\overline{b}\overline{a})^{-2}
$%
.

Dans les cas
${\romannumeral 3)\/}$
et
${\romannumeral 4)\/}$
les choix sont deux \`a deux non conjugu\'es%
\footnote{car 
$S, S^{-1}, BA\!=\!M,  M^{-1}, M^{2}, M^{-2}$
sont distinctes dans l'ab\'elianis\'e de
$B_{3}/Y=SL(2, {\Bbb Z})$%
.
}.

\finc
\vskip3mm
\Dem
Si
$t_{n}(w)=(k, m)\geq (2, 1)$
alors
$w\!=\!sWc^{n},\ \hbox{\rm ou\/}\ sWsc^{n}$
et le conjugu\'e
$w^{s}\!=\!Wsc^{n}$
ou
$w^{s}\!=\!Wc^{n+1}$
est de taille normale
$(k-1,m)$
ou
$(k-3, m)$%
.

Si
$w=Wsc^{n}$
a une \'ecriture de taille
$(1, m),\, m\!>\!1$
alors
$W$
a l'une des formes~:
$$W\!=\!a^{-1}W'a^{-1},\, a^{-1}W'b,\, bW'a^{-1},\, bW'b$$
avec
$l(W')\!=\!l(W)-2\!<\!l(W)$
et dans les quatre cas les conjugu\'es~:
$$w^{a^{-1}}\!=\!aa^{-1}W'a^{-1}abac^{n}a^{-1}=W'bc^{n},\
w^{a^{-1}}\!=\! aa^{-1}W'babac^{n}a^{-1}=W'sc^{n}$$
$$w^{b}\!=\!b^{-1}bW'a^{-1}abac^{n}b=W'sc^{n},\
w^{b}=b^{-1}bW'babac^{n}b=W'a^{-1}c^{n+1}$$
ont une \'ecriture de taille normale inf\'erieure~:
$$(0, m-1),\ (1, m-2),\ (1, m-2),\ (0, m-1)\ .$$

Dans le cas
${\romannumeral 4)\/}$
si
$n$
est pair~:
$a^{-1}s^{2n+1}\!\equiv\!a^{-1}s\!=\!a^{-1}aba\!=\!ba$
mod
$Y$%
,
$bs^{2n+1}\!\equiv\!bs\!=\!baba\!=\!(ba)^{2}$
mod
$Y$
et si
$n$
est impair
$a^{-1}s^{2n+1}\!\equiv\!a^{-1}s^{3}\!=\!a^{-1}ababababa\!=\!babababa\!=\!(ba)^{4}$%
mod
$K$%
,
$bs^{2n+1}\!\equiv\!bs^{3}\!=\!bababababa\!=\!(ba)^{5}$
mod
$K$%
.\hfill\break
\null\hfill\findem

\Thc Formules dans
$Aut^{+}(F_{2}) = Int(F_{2})\rtimes_{\psi}B_{3})$|
${\romannumeral 1)\/}\quad
\sigma\varphi_{x}=\varphi_{\sigma(x)}\sigma$%

${\romannumeral 2)\/}\quad (\varphi_{x}\sigma)(\varphi_{y}\tau)=
\varphi_{x\sigma(y)}\sigma\tau$%

${\romannumeral 3)\/}\quad (\varphi_{x}\sigma)^{\varphi_{y}\tau}=
\varphi_{\tau^{-1}(y^{-1}x\sigma(y))}\sigma^{\tau}$
.
\finc
\vskip2mm
\Dem
${\romannumeral 2)\/}$
suit de
${\romannumeral 1)\/}$
et
${\romannumeral 3)\/}$
suit  de
${\romannumeral 2)\/}$
et
${\romannumeral 1)\/}$%
~:
$$
\varphi_{x}\sigma\varphi_{y}\tau=\varphi_{x}\varphi_{\sigma(y)}\sigma\tau=
\varphi_{x\sigma(y)}\sigma\tau
$$
$$
(\varphi_{y}\tau)^{-1}\!\varphi_{x}\sigma\varphi_{y}\tau\!=\!%
\tau^{-1}\!\varphi_{y^{-1}}\varphi_{x\sigma(y)}\sigma\tau\!=\!%
\tau^{-1}\!\varphi_{y^{-1}x\sigma(y)}\sigma\tau\!=\!
\varphi_{\tau^{-1}\!(y^{-1}x\sigma(y))}\tau^{-1}\sigma\tau
$$
Quant \`a
${\romannumeral 1)\/}$%
, pour tout
$t\!\in\!F_{2}$
on a~:
$$\sigma\varphi_{x}(t)\!=\!\sigma(xtx^{-1})\!=\!%
\sigma(x)\sigma(t)\sigma(x)^{-1}\!=\!\varphi_{\sigma(x)}(\sigma(t))\!=\!%
(\varphi_{\sigma(x)}\sigma)(t)\ .$$
\hfill\findem

\Thc Remarque|
Les \'el\'ements suivants
$\sigma,\, \iota,\, \zeta$
de
$Aut^{+}(F_{2})$%
~:
$$
\sigma=\varphi_{vu}s^{2}:\, (u, v)\mapsto (u^{-1}, v^{-1})
$$
$$
\iota=\varphi_{u}s:\, (u, v)\mapsto (v, u^{-1}),\
(\iota^{2}=\sigma)
$$
$$
\zeta:=\varphi_{u}\beta s: (u, v)\mapsto(v, v^{-1}u^{-1}),\quad
\zeta^{-1}(u, v)=(v^{-1}u^{-1}, u)
$$
sont d'ordre fini, respectivement
$2, 4$
 et
 $3$
 et d'images dans
$SL(2, {\Bbb Z\/})$%
~:
$$
\rho(\sigma)\!=\!-I\!=\!\begin{pmatrix}-1& 0\cr 0&\!\!-1\cr\end{pmatrix},
\rho(\iota)\!=\!S\!=\!\begin{pmatrix}0&\!-1\cr 1& 0\cr\end{pmatrix}, 
\rho(\zeta)\!=\!\begin{pmatrix}0& \!-1\cr 1&\!-1\cr\end{pmatrix},
\rho(\zeta^{-1}\!)\!=\!\begin{pmatrix}-1& 1\cr-1& 0\cr\end{pmatrix}.
$$

De plus
$\zeta^{\sigma}=\varphi_{v}\zeta$
et
$\sigma^{\zeta^{-1}}\!=\varphi_{v^{-1}}\sigma$%
.
\finc
\vskip2mm
\Dem
D'apr\`es la Proposition de l'Appendice B  on a les relations~:
$s^{2}(u, v)\!=\!(u^{-1}\!\!\,, v^{-1}\!)^{vu}\!\!,\ s(u, v)\!=\!(v, u^{-1}\!)^{u}$%
, d'o\`u les cas de
$\sigma$
et
$\iota$
et~:\hfill\break
%
$
\zeta(u, v)\!=\!\varphi_{u}\beta((v, u^{-1}\!)^{u})\!=\!%
\varphi_{u}(v^{vu}, (u^{-1}\!v^{-1}\!)^{vu})\!=\!(v^{v}, (u^{-1}\!v^{-1}\!)^{v})\!=\!%
(v, v^{-1}\!u^{-1}\!)
$%
,\hfill\break  expression en
$u$
et
$v$
de la permutation cyclique des g\'en\'erateurs
de la pr\'e\-sen\-tation sym\'etrique
$
F_{2}\!=\langle u, v, w\,|\, uvw\rangle
$
d'o\`u l'ordre de
$\zeta$
et l'expression de
$\zeta^{-1}$%
.\hfill\carre

$$
\zeta^{\sigma}(u, v)\!=\!\sigma^{-1}\zeta\sigma(u, v)\!=\!%
\sigma\zeta(u^{-1}\!, v^{-1}\!)\!=\!%
\sigma(v^{-1}\!, uv)\!=\!(v, u^{-1}\!v^{-1}\!)\!=\!%
\varphi_{v}\zeta(u, v) 
$$
$$
\sigma^{\zeta^{-1}}\!(u, v)\!=\!\zeta\sigma\zeta^{-1}(u, v)\!=\!%
\zeta\sigma(v^{-1}u^{-1}, u)\!=\!\zeta(vu, u^{-1})\!=\!(v^{-1}\!u^{-1}\!v, v^{-1})\!=\!%
\varphi_{v^{-1}}\sigma(v)\ .
$$
\null\hfill\findem

\Thc coCorollaire|
${\romannumeral 1)\/}$
Les \'el\'ements de torsion  de
$SL(2, {\Bbb Z})$
sont d'ordre
$6, 4, 3, 2$%
, une classe de conjugaison pour ceux d'ordre
$2$%
, deux pour les autres~:
$$
M\!=\!BA\!=\!\!\begin{pmatrix}1&\!-1\cr 1& 0\cr\end{pmatrix}\!,
M^{-1}\!=\!\!\begin{pmatrix}0& 1\cr-1&\!\!1\cr\end{pmatrix} 
S\!, S^{-1}\!,
M^{2}\!=\!\!\begin{pmatrix}0& -1\cr 1& -1\cr\end{pmatrix}\!,
M^{-2}\!=\!\!\begin{pmatrix}-1& 1\cr-1& 0\cr\end{pmatrix}\!,
-I\ .
$$

${\romannumeral 2)\/}$
Tout relev\'e dans
$Aut^{+}(F_{2})$%
, avec m\^eme ordre, d'un \'el\'ement d'ordre
$2,\, 4$
ou
$3$
de
$SL(2, {\Bbb Z\/})$
est respectivement conjugu\'e \`a
$\sigma$%
,
$\iota$
ou
$\iota^{-1}$%
,
$\zeta$
ou
$\zeta^{-1}$%
.

Il n'y a pas d'\'el\'ement d'ordre 
$6$
dans
$Aut^{+}(F_{2})$%
.
\finc
\vskip2mm
\Dem
${\romannumeral 1)\/}$
est une paraphrase du {\petcap Compl\'ement\/}.\hfill\carre

${\romannumeral 2)\/}$
Si
$\gamma\in Aut^{+}(F_{2})$
est d'ordre
$k$
et
$\gamma'\in Aut^{+}(F_{2})$
est tel que
$\rho(\gamma')=\rho(\gamma)$
alors, d'apr\`es Nielsen 
${\romannumeral 3)\/}$
il y a
$x\in F_{2}$
tel que
$\gamma'=\varphi_{x}\gamma$%
.

D'apr\`es la
{\petcap Formule\/}
${\romannumeral 2)\/}$%
,
$\gamma'$
est d'ordre
$k$
si et seulement si~:
$$
x\in X_{\gamma}:=\{x\in F_{2}\, |\,
x\gamma(x)\cdots\gamma^{k-2}(x)\gamma^{k-1}(x)=1\}
$$
et, d'apr\`es la
{\petcap Formule\/}
${\romannumeral 3)\/}$
et Nielsen
${\romannumeral 3)\/}$%
, l'\'el\'ement
$\varphi_{y}\tau$
conjugue
$\gamma''=\varphi_{x'}\gamma$
\`a
$\gamma'$
si et seulement si~:
$\gamma^{\tau}=\varphi_{z}\gamma$
pour un
$z\in F_{2}$
et
$\tau^{-1}(y^{-1}x'\gamma(y))z=x$%
.

Si
$\tau\!=\!1$
on a
$z\!=\!1$
et
$x$
est dans l'orbite de
$x'$
pour l'action \`a droite~:
$$\cdot_{\gamma} : X_{\gamma}\times F_{2} \rightarrow X_{\gamma},\
x\cdot_{\gamma} w=w^{-1}x\gamma(w)\ .$$

D'apr\`es
${\romannumeral 1)\/}$
et Nielsen
${\romannumeral 3)\/}$%
, les \'el\'ements d'ordre
$2,\, 3,\, 4$
de
$Aut^{+}(F_{2})$
sont conjug\'es \`a l'un de~:
$
\varphi_{x_{\sigma}}\sigma,\ (\varphi_{x_{\zeta}}\zeta)^{\pm1},\ (\varphi_{x_{\iota}}\iota)^{\pm1}
$
avec des
$x_{\gamma}\in X_{\gamma},\ \gamma=\sigma,\, \zeta,\, \iota$%
.

Soit
$x_{\gamma}\!\in\!X_{\gamma}$
de longueur 
$l(x_{\gamma})\!=\!\min\{l(w\cdot_{\gamma}\!x_{\iota}\, |\, w\!\in\!F_{2}\}\!=\!l$
minimale dans l'orbite de
$x_{\gamma}$
pour l'action de
$\cdot_{\gamma}$
avec
$\gamma\!=\!\sigma, \zeta,\, \iota$%
.
\vskip2mm
\Demdsp {es cas d'ordre non divisible par
$3$
dans
${\romannumeral 2)\/}$%
}|

Si
$l=0,\ \varphi_{x_{\gamma}}\gamma$
est, suivant le cas,
$\sigma$
ou
$\iota^{\pm1}$%
.\hfill\carre

Sinon,
$\gamma$
envoyant les g\'en\'erateurs sur des \'el\'ements de longueur
$1$,\ $\gamma(x_{\gamma})$%
, en mot r\'eduit, s'obtient en rempla\c cant dans
$x_{\gamma}$
chaque lettre par son image par
$\gamma$.

Le mot
$x_{\gamma}\gamma(x_{\gamma})\cdots\gamma^{k-1}(x_{\gamma})=1$
n'\'etant pas r\'eduit, soit
$x_{\gamma}$
est de longueur
$l(x_{\gamma})\geq 2$
avec
$a, z$
pour premi\`ere et derni\`ere lettre v\'erifiant
$\gamma(a)=z^{-1}$%
, donc
$l(x_{\gamma}\cdot_{\gamma} a)=l(x_{\gamma})-2$%
, contredisant la minimalit\'e de
$l(x_{\gamma})$
dans l'orbite de
$x_{\gamma}$%
.\hfill\carre

Donc
$l(x_{\gamma})=1,\ x_{\gamma}=a$
et
$\gamma(a)=a^{-1}$%
, impossible si
$\gamma=\iota$%
, et si
$\gamma=\sigma$%
, comme
$\iota$
agit transitivement sur
$\{u^{\pm1}, v^{\pm1}\}$
les quatre cas
$x_{\sigma}=\varphi_{t}
\sigma,\, t\in\{u^{\pm1}, v^{\pm1}\}$
sont conjugu\'es entre eux par les puissances de
$\iota: (\varphi_{t}\sigma)^{\iota^{k}}\!=\!\varphi_{\iota^{-k}(t)}\sigma,\,  k\!=\!1, 2, 3, 4$
et, comme%
\footnote{par la derni\`ere relation de la {\bf Remarque\/}.
}\
$\sigma^{\zeta^{-1}}=\varphi_{v^{-1}}\sigma$%
, l'\'el\'ement
$\varphi_{x_{\sigma}}\sigma$
est aussi conjugu\'e \`a
$\sigma$%
. \hfill\carre\findem
\vskip2mm

\Demdsp{u cas d'ordre 
$3$
dans
${\romannumeral 2)\/}$%
}|

Par confort notationnel
on identifie les g\'en\'e\-rateurs
de la pr\'esentation sym\'etrique de
$F_{2}\!=\!\langle u, v, w\, |\, uvw\rangle$
\`a
$\overline{1}, \overline{2}, \overline{3}\in{\Bbb Z\/}/3{\Bbb Z\/}$%
, ainsi~:
$$
F_{2}=\langle{\Bbb Z\/}/3{\Bbb Z\/}\, |\, (i-1)i(i+1)\rangle,\
\zeta: F_{2}\rightarrow F_{2}, \zeta(i)=i+1\ .
$$

Comme pour tout g\'en\'erateur
$i\in{\Bbb Z\/}/3{\Bbb Z\/}$%
, on a~:
$
i^{-1}=(i+1)(i-1)
$%
, tout \'el\'ement de 
$w\!\in\!F_{2}$
a une \'ecriture
$w\!=\!i_{1}\cdots i_{n}$
en uniquement les g\'en\'erateurs (et pas leurs inverses)
$\overline{1}, \overline{2}, \overline{3}\in{\Bbb Z\/}/3{\Bbb Z\/}$%
, \'ecriture  unique si elle est {\it r\'eduite\/}, c.\`a.d. ne contient,
pour un
$i\!\in\!{\Bbb Z\/}/3{\Bbb Z\/}$%
, le sous-mot 
$
(i-1)i(i+1)
$%
, on note alors
$l^{+}(w)\!=\!n$%
, sa {\it longueur\/}.

Soit
$x\!=\!i_{1}\cdots i_{n}\!\in\!X_{\zeta}$
avec
$l^{+}(x)$
minimale dans l'orbite de
$x$
pour l'action
$\cdot_{\zeta}$

\Thc Affirmation|
$n=l^{+}(x)\leq1$%
\finc
\Dem
On a
$x\zeta(x)\zeta^{2}(x)\!=\!1$
donc, soit
$l^{+}(x)\!\leq\!1$%
, soit il y a  r\'eduction entre%
\footnote{l'image par
$\zeta^{-1}$
d'une r\'eduction entre
$\zeta(x)$
et
$\zeta^{2}(x)$
est une r\'eduction entre
$x$
et
$\zeta(x)$%
.
}\
$x$
et
$\zeta(x)$%
.
Comme
$(i_{1}i_{2})(i_{1}+1)(i_{2}+1)$
ne se r\'eduit pas, on a
$l(x)\!=\!n>2$%
.

Comme
$\zeta(x)\!=\!(i_{1}+1)\cdots (i_{n}+1),\ i_{n}\!=\!i_{1}$
et soit
$i_{2}\!=\!i_{1}+1$%
, soit
$i_{n-1}\!=\!i_{1}-1$%
.

Dans le premier cas%
\footnote{Si
$n\!=\!3, 4, 5,\, i_{3}\cdots i_{n-1}$%
(resp.
$i_{2}\cdots i_{n-2}$%
)
sont
$\emptyset, i_{3}, i_{3}i_{n-1}\!=\!i_{3}i_{4}\ (\hbox{\rm resp.\/}\,
\emptyset, i_{2}, i_{2}i_{n-2}\!=\!i_{2}i_{3})$%
.
}\
$ x\cdot_{\zeta}(i_{n}-1)^{-1}\!=\!%
i_{1}(i_{1}+1)i_{3}\cdots i_{n-1}i_{1}\cdot_{\zeta}(i_{1}-1)^{-1}\!=\!%
(i_{1}-1)i_{1}(i_{1}+1)i_{3}\cdots i_{n-1}\!=\!i_{3}\cdots i_{n-1}$%
, contredisant
$l^{+}(x)$
minimal dans l'orbite.

Dans le second 
$x\cdot_{\zeta}i_{1}\!=\!i_{1}i_{2}\cdots i_{n-2}(i_{1}-1)i_{1}\cdot_{\zeta}i_{1}\!=\!%
i_{2}\cdots i_{n-2}(i_{1}-1)i_{1}(i_{1}+1)\!=\!i_{2}\cdots i_{n-2}$%
, contredisant encore
$l^{+}(x)$
minimal dans l'orbite de
$x$%
.\hfill\findem

Ainsi soit
$x\!=\!1$
et
$\zeta'\!=\!\zeta$%
, soit
$\zeta'\!=\!\varphi_{\overline{i}}\zeta$%
, cependant%
\footnote{d'apr\`es la {\bf Remarque.\/}
}\
$\zeta^{\sigma}\!\!=\!\varphi_{\overline{1}}\zeta$
et, pour
$k\!=\!1, 2, 3,\, (\varphi_{\overline{i}}\zeta)^{\zeta^{-k}}\!=\!%
\varphi_{\zeta^{k}(\overline{i})}\zeta\!=\!%
\varphi_{\overline{i+k}}\,\zeta$
et
$\zeta'$
est encore conjugu\'e \`a
$\zeta$%
.\hfill\qed\qed
\vskip2mm
\Demdsp{e la non existence d'\'el\'ement d'ordre
$6$
dans
${\romannumeral 2)\/}$%
}|

Soit
$\omega$
d'ordre
$6$
de
$Aut^{+}(F_{2})$%
. Les \'el\'ements d'ordre
$3$
de
$Aut^{+}(F_{2})$
\'etant conjugu\'es \`a
$\zeta$
ou
$\zeta^{-1}$
et, dans
$SL(2, {\Bbb Z\/})$%
, ceux d'ordre
$6$
 \`a
$M$
ou
$M^{-1}$%
, on suppose 
$\omega$
d'ab\'elianis\'e 
$M$
et
de carr\'e
$\omega^{2}\!=\!\zeta\!=\!\varphi_{u}\beta s$%
, le cas
$M^{-1}$
suit par passage \`a l'inverse.

Ainsi
$\omega=\varphi_{x}\beta\alpha$
et~:
$$\omega^{2}=\varphi_{x}\beta\alpha\varphi_{x}\beta\alpha\!=\!%
\varphi_{x}\varphi_{\beta\alpha(x)}\beta\alpha\beta\alpha\!=\!%
\varphi_{x\beta\alpha(x)}\beta s\!=\!%
\varphi_{x\beta\alpha(x)}\varphi_{u^{-1}}\zeta\!=\!%
\varphi_{x\beta\alpha(x)u^{-1}}\zeta\ .
$$

Mais
$N\!=\!I+M\!=\!\begin{pmatrix}2& \!-1\cr 1& 1\cr\end{pmatrix}
$
ne peut avoir
$\overline{u}\!=\!e_{1}\!=\!\begin{pmatrix}1\cr 0\cr\end{pmatrix}$
dans son image%
\footnote{l'unique pr\'eimage par
$N$
de
$e_{1}$
dans
${\bf Q\/}^{2}$
est
${{1}\over{3}}e_{1}-{{1}\over{3}}e_{2}
\not\in{\bf Z\/}^{2}$%
. 
}, \`a plus forte raison, dans le groupe libre
$F_{2}=\langle u, v\rangle$%
, l'\'equation~:
$x\beta\alpha(x)u^{-1}\!=\!1$%
, dont
$(I+M)(X)-U\!=\!0$
est l'\'equation induite dans l'ab\'elianis\'e
${\Bbb Z\/}^{2}$%
, n'a pas de solution, donc il n'y a pas d'hypoth\'etique
$\omega$
d'ordre
$6$
dans
$Aut^{+}\!(F_{2})$%
.\qed\qed\qed


\section*{APPENDICE D\\
Le groupe des tresses di\'edrales \`a trois brins
$DB_{3}$
et
$Aut(F_{2})$%
.
}

Le goupe \`a deux \'el\'ements
$\Delta=\{1, d\}$
agit sur le groupe des tresses \`a trois brins
$B_{3}$
par
$a^{d}=b^{-1}, b^{d}=a^{-1}$%
.
\Defn
Le groupe des {\it tresses di\'edrales \`a trois brins\/}
$DB_{3}$
est l'extension triviale de
$\Delta$
par
$B_{3}$
muni de l'action pr\'ec\'edente.~:
$$
1\rightarrow B_{3}\rightarrow DB_{3}=B_{3}\rtimes\Delta\rightarrow \Delta\rightarrow 1\ ,
$$
il a donc la pr\'esentation~:
$$
DB_{3}=\langle a, b, d\,|\, d^{2},  aba(bab)^{-1}, bdad\rangle\ .
$$

\Thc Proposition|
L'endomorphisme
$\delta: F_{2}\rightarrow F_{2}$
d\'efini par
$\delta(u, v)=(v, u)$
est un isomorphisme d'ordre
$2$
v\'erifiant les relations 
$$
\alpha^{\delta}=\beta^{-1},\quad \beta^{\delta}=\alpha^{-1}\ ,
$$
il permet donc d'\'etendre le morphisme~:
$$
\Psi: B_{3}\rightarrow Aut^{+}(F_{2}),\ \Psi(a)=\alpha,\ \Psi(b)=\beta
$$
de l'appendice B en un morphisme injectif~:
$$
\Phi: DB_{3}\rightarrow Aut(F_{2}), \ \Phi(a)=\Psi(a)=\alpha,\
\Phi(b)=\Psi(b)=\beta,\ \Phi(d)=\delta\ .
$$
\finc

\Dem
Comme
$\delta^{2}=1$%
, il suffit d'\'etablir
$\alpha^{\delta}=\beta^{-1}$%
~:\hfill\break
$
\alpha^{\delta}(u, v)=\delta^{-1}\alpha\delta(u, v)=\delta\alpha(v, u)=
\delta(u^{-1}v, u)=(v^{-1}u, v)=\beta^{-1}(u, v)
$%
.\hfill\findem

Il sera alors ais\'e \`a la lectrice d'\'etablir, pour 
$Aut(F_{2})$%
, la version originelle~:

\Thc Th\'eor\`eme$^{\prime}$
{\rm (de Nielsen, 1917)}|
${\romannumeral 1})$
Un endomorphisme surjectif du groupe libre \`a deux g\'en\'erateurs
$F_{2}$
est un automorphisme.

${\romannumeral 2)}$
Le groupe des automorphismes de
$F_{2}$
est produit semi-direct amalgam\'e~:
$$Aut(F_{2}) = Int(F_{2})\rtimes_{<\varphi_{[u, v]}\,=\,\Phi(s^{4})>}\Phi(DB_{3})\ .$$

${\romannumeral 3)}$
Le groupe
$Int(F_{2})$
des automorphismes int\'erieurs de
$F_{2}$
est le noyau de~:
$$\rho: Aut(F_{2})\rightarrow Aut({\Bbb Z}^{2})=GL(2, {\Bbb Z})\ .$$

${\romannumeral 4)}$
Un endomorphisme
$\varphi$
du groupe libre \`a deux g\'en\'erateurs
$F_{2}$
d\'efini par
$\varphi(t_{i})\!=\!X_{i}, i\!=\!1, 2$
est un automorphisme  de
$F_{2}$
si et seulement si il pr\'eserve \`a conjugaison et inverse pr\`es
le commutateur
des g\'en\'erateurs~:
$$[X_{1}, X_{2}]=X_{1}^{-1}X_{2}^{-1}X_{1}X_{2}=([t_{1}, t_{2}]^{\pm1})^{w}=
w^{-1}[t_{1}, t_{2}]^{\pm1}w\ .
\leqno{(C'w)}$$
De plus, il est dans
$\Phi(DB_{3})$
ssi
$[X_{1}, X_{2}]\!=\![t_{1} ,t_{2}]^{\pm1}$
, c. \`a d. on a 
$(C'1)$
\finc

\noindent
et la d\'etermination \`a conjugaison pr\`es des torsions de
$GL(2, {\Bbb Z\/})$
et
$Aut(F_{2})$%
~:

\Thc coCorollaire|
${\romannumeral 1)\/}$
Les \'el\'ements de torsion  de
$GL(2, {\Bbb Z})$
sont d'ordre
$2, 3, 4$%
\hfill\break
ou
$6$%
, trois classes de conjugaison pour ceux d'ordre
$2$%
, une seule pour les autres~:
$$-I\!,D\!=\!\begin{pmatrix}0& 1\cr 1& 0\cr\end{pmatrix}, SD\!=\!%
\begin{pmatrix}-1& 0\cr 0&1\cr\end{pmatrix},\,
M^{2}\!=\!\!\begin{pmatrix}0& -1\cr 1& -1\cr\end{pmatrix},\, S\
\hbox{\rm ou\/}\
M\!=\!\!\begin{pmatrix}1&\!-1\cr 1& 0\cr\end{pmatrix}\ .
$$

${\romannumeral 2)\/}$
Tout relev\'e dans
$Aut^{+}(F_{2})$%
, avec m\^eme ordre, d'un \'el\'ement d'ordre
$2, 3, 4$
de
$GL(2, {\Bbb Z\/})$
est respectivement conjugu\'e \`a
$\sigma,\ \delta,\, \iota\delta,\, \zeta$
ou 
$\iota$%
.

Il n'y a pas d'\'el\'ement d'ordre
$6$
dans
$Aut(F_{2})$%
.\hfill\qed\findem 
\finc


\section*{APPENDICE E\\
Un exercice de th\'eorie des groupes~: produits semi-directs amalgam\'es.}

\parc
On rapelle que, si
$K, L$
sont des groupes et
$\psi: L\!\rightarrow\!Aut(K)$
est un morphisme de groupes, le produit semi-direct de
$K$
par
$L$
suivant
$\psi$
est
$K\!\rtimes_{\psi}\!L$%
, le produit
$K\!\times\!L$
muni de la loi~:
$$(k, l)\cdot(k', l')\!=\!(k\psi(l)(k'), ll')\ .$$

Il contient
$K\!\times\!\{1\}, \{1\}\!\times\!L\subset K\!\rtimes_{\psi}\!L$
comme sous-groupes qu'abusivement, on note encore
$K\!=\!K\!\times\!\{1\}, L\!=\!\{1\}\!\times\!L< K\!\rtimes_{\psi}\!L$
avec
$K\cap L\!=\!\{1\}$%
, le premier
$K$
\'etant normal et la restriction
$\pi_{|L}: L\subset K\!\rtimes_{\psi}\!L\rightarrow K\!\rtimes_{\psi}\!L/K$
\`a
$L$
du morphisme quotient est un isomorphisme.

De plus
$(k, l)^{-1}\!=\!(\psi(l^{-1})(k^{-1}), l^{-1})$
et, si
$l\in\!L,$
\hfill\break
$\psi(l): K\!=\!K\!\times\!\{1\}\!\rightarrow\!K\!\times\!\{1\}\!=\!K$
est la restriction \`a
$K\!\times\!\{1\}$
de la conjugaison par
$l$
dans
$K\!\rtimes_{\psi}\!L, (k, 1)\mapsto (1, l)(k, 1)(1, l)^{-1}\!=\!\psi(l)(k)$%
.
\finc

Cet appendice donne  une version, sans l'hypoth\`ese
$K\cap L\!=\!\{1\}$%
, de la classique~:
\Thc Proposition|
{\sl Soit
$K< G$
un sous-groupe  d'un groupe
$G$
et
$L<N(K)$
un  sous-groupe du normalisateur de
$K$
avec
$K\cap L\!=\!\{1\}$%
. On note
$$\psi: L\!\rightarrow\!Aut(K), \psi(l)(k)\!=\!lkl^{-1}$$

Alors 
$K\!\rtimes_{\psi}\!L$
est, par
$(k, l)\mapsto kl$%
, isomorphe au sous-groupe
$KL\!\subset\!G$%
.}
\finc

Soit
$K, L, M$
trois groupes,
$\psi: L\!\rightarrow\!Aut(K), \kappa: M\!\rightarrow\!K,
\lambda:M\!\rightarrow\!L$
trois morphismes de groupes avec
$\kappa, \lambda$
injectifs,
$\lambda(M)\triangleleft L$
normal dans
$L$
et
$\psi$
respectant l'image
$\kappa(M)$%
, c.\`a.d. induisant
$\psi_{M}: L\!\rightarrow\!Aut(\kappa(M))$%
, et tant
$\psi_{M}$
que
$\psi\circ\lambda$
sont int\'erieurs au sens fort suivant%
\footnote{notant abusivement 
$\lambda^{-1}: \lambda(M)\rightarrow M$
l'inverse de la bijection
$\lambda': M\rightarrow \lambda(M), m\mapsto \lambda(m)$%
.
}~:
$$\psi_{M}(l)(\kappa(m))\!=\!\kappa\circ\lambda^{-1}(l\lambda(m)l^{-1}),\
\psi(\lambda(m))(k)\!=\!\kappa(m)k\kappa(m)^{-1}\ .$$

On note
$\mu: M\rightarrow K\!\rtimes_{\psi}\!L, \mu(m)\!=\!%
(\kappa(m)^{-1}, \lambda(m))$%
.

Identifiant, au moyen de
$\kappa$
et
$\lambda$%
, le groupe
$M$
\`a un {\it sous-groupe commun\/} \`a
$K$
{\it et\/}
$L$%
, ainsi que
$K, L$
aux sous-groupes
$K\!\times\!\{1\}, \{1\}\!\times\!L$
de
$K\!\rtimes_{\psi}\!L$%
, on a~:
$$\mu: M\!\rightarrow\!K\!\rtimes_{\psi}\!L,\ \mu(m)\!=\!(m^{-1}, m)\ ,$$
$$\psi_{M}(l)(m)\!=\!lml^{-1}\!=:\!^{l}m,\quad
\psi(m)(k)\!=\!mkm^{-1}\!=:\!^{m}k\ ,$$
$$(k, l)\!=\!kl,\quad (k, l)(k', l')\!=\!(k(^{l}k'), ll'),\quad
(k, l)^{-1}\!=\!(l^{-1}k^{-1}l, l^{-1})\!=:\!((k^{-1})^{l}, l^{-1})\ .$$

\Thc Lemme|
$\mu$
est un morphisme  d'image
$\mu(M)\triangleleft K\!\rtimes_{\psi}\!L$%
, un sous-groupe normal.
\finc
\Defn
Le groupe quotient
$K\!\rtimes_{\psi}\!L/\mu(M)$%
, not\'e
$K\!\rtimes_{\psi, \eta=\kappa}\!L$
(ou%
\footnote{comme ci-dessus, identifiant
$M$
\`a
$\kappa(M),\cdots$
}\
$K\!\rtimes_{\psi, M}\!L$%
),
est le {\it produit semi-direct suivant\/}
$\psi$
{\it de\/}
$K$
{\it et\/}
$L$
{\it amalgam\'e par\/}
$\eta\!=\!\kappa$
{\it (ou sur\/}
$M$%
{\it ).}
\Thc Corollaire|
Soit
$K< G$
un sous-groupe  d'un groupe
$G$
et
$L<N(K)$
un  sous-groupe du normalisateur de
$K$%
. Alors, si
$\psi: L\rightarrow Aut(K), \psi(l)(k)\!=\!lkl^{-1}$
est la restriction \`a
$L<G$
de la conjugaison dans
$G$%
, l'application~:
$$p: K\!\rtimes_{\psi}\!L\rightarrow KL<G,\ (k, l)\mapsto kl$$
est un morphisme de noyau
$M\!=\!K\!\cap\!L$
et induit un isomorphisme~:
$$\overline{p}: K\!\rtimes_{\psi, M}\!L\rightarrow KL\ .$$

\finc
\parc
\Dem
on a 
$p(k, l)p(k', l')\!=\!klk'l'\!=\!klk'l^{-1}ll'\!=\!p((k, l)(k', l'))$%
, ainsi
$p$
est  morphisme de groupes.\hfill\carre\break
$1\!=\!p(k, l)\!=\!kl$
ssi
$K\!\ni\!k\!=\!l^{-1}\!\in\!L$%
, donc
$l\in K\!\cap\!L\!=\!M$%
, c.\`a.d.
$(k, l)\!\in\!\mu(L)$%
.\hfill\carre\findem
\vskip2mm
\Demd u Lemme|%
$\mu(m)\mu(m')\!=\!(m^{-1}, m)(m'^{-1}, m')\!=\!(m^{-1}(^{m}m'^{-1}), mm')\!=\!$
$(m'^{-1}m^{-1}, mm')\!=\!((mm')^{-1}, mm')\!=\!\mu(mm')$
donc
$\mu$
est un morphisme.\hfill\carre

$$^{kl}\mu(m)=kl(m^{-1}, m)(kl)^{-1} = (k(^{l}m^{-1}), lm)l^{-1}k^{-1} =
(k(^{l}m^{-1}), lm)(l^{-1}k^{-1}ll^{-1}) =$$
$$=(k(^{l}m^{-1})(^{lm}((k^{-1})^{l})), lml^{-1}) =
(klm^{-1}l^{-1}lml^{-1}k^{-1}lm^{-1}l^{-1}, ^{l}m) =$$
$$(^{l}m^{-1}, ^{l}m) = \mu(^{l}m)\in\mu(M)\ .$$
et le sous-groupe
$\mu(M)$
est normal dans
$K\!\rtimes_{\psi}\!L$%
.\hfill\carre\findem

\finc


\section*{{\petcap Commentaires bibliographiques}}

Les seuls calculs
$ABA\!=\!BAB, f_{n+1}f_{n-1}\!=\!f_{n}$
et
$\pi_{0}(A)\!=\!a, \pi_{0}(B)\!=\!b$
plut\^ot que ceux et les  r\'ecurrences \`a base d'algorithme d'Euclide
du traitement classique (p. 9-19 de \cite{Ra})
profitent de ce que  la suffisance des relations est donn\'ee par
la pr\'esentation du groupe des 
classes d'isotopie de diff\'eomorphismes de la sph\`ere plate
$({\Bbb S}, \{m, n, p, \infty\}, \{\infty\})=
({\Bbb T}/\{\pm I\}, {\Bbb T}_{2}, \{0\})$
respectant points d'ordre
$2$
et
$1$%
, pr\'esentation d\'eduite de celle du
groupe des tresses \`a trois brins  d'Artin (\cite{A}).

Plus de vingt ans apr\`es%
\footnote{{\sl Theory of Braids} Ann. of Math. 48 (1947), p. 101-126.%
}%
$\!$%
, Artin conc\'edait que (dans son article de 1925)
``{\sl Most of the proofs are entirely intuitive}'',
renvoyant pour les relations de
$B_{n}$
au quasiment simultan\'e%
\footnote{{\sl The algebraic braid group.} Ann. Math. 48 (1947), p. 127-136.
}\
traitement alg\'ebrique de F. Bohnenblust via l'action de
$B_{n}$
sur le groupe libre
$F_{n}$
(le groupe fondamental du plan priv\'e de
$n$
points).

Depuis, l'\'etablissement de ces relations%
\footnote{tant celles du groupe des tresses que le passage au groupe des classes d'isotopies
d'hom\'eo\-mor\-phismes ou diff\'eomorphismes respectant une partie finie
de la sph\`ere.
}\
a fait couler beaucoup d'encre.

Un traitement clair et rigoureux est donn\'e par l'\'ecole d'Orsay
via la m\'ethode des espaces fonctionnels%
\footnote{L'espace des plongements de
$n$
points dans le plan est stratifi\'e avec une seule strate ouverte
(les
$n$
points d'abscisses deux \`a deux distintes) qui est contractile,
les g\'en\'erateurs du groupe fondamental correspondant aux travers\'ees des strates de codimension
$1$
et  les  relations aux contour\-nements de celles de codimension
$2$
(des strates de  co\"{\i}ncidences d'abscisses).
}\
et les th\'eor\`emes de fibration de Cerf (\cite{Ce}).%
\footnote{voir la r\'edaction de Valentin Po\'enaru, Expos\'e 2, p. 21-31 de \cite{FLP}.
}

Dans son introduction \`a la K-th\'eorie alg\'ebrique (\cite{Mi}) Milnor
prouve%
\footnote{chaque \'etape est simple, une fois assimil\'e le d\'ebut
du \S5 et beaucoup de \S8, \S9 et \S10\dots!
}\
que le noyau du morphisme canonique
$\phi\!:\!St(2, {\Bbb Z})\!\rightarrow\!\!SL(2, {\Bbb Z})$
du groupe de Steinberg%
\footnote{$St(2, {\bf Z})=<x_{12}, x_{21}\,|\, x_{21}^{x_{12}^{-1}x_{21}x_{12}^{-1}}=x_{12}^{-1}>$
et
$\phi(x_{12})=e_{12}=A^{-1}, \phi(x_{21})=e_{21}=B$%
.
}\ 
sur le groupe sp\'ecial lin\'eaire en dimension
$2$
est un groupe cyclique infini.

Exprimant la pr\'esentation de
$St(2, {\Bbb Z})$
au moyen de
$\alpha=x_{1, 2}^{-1}, \beta=x_{2, 1}$%
, il recon\-na\^{\i}t
$St(2, {\Bbb Z})\!\simeq\!B_{3}$%
, le noyau de
$\phi$
\'etant le sous-groupe d'indice
$2$
du centre de ce groupe, non  de tresses,
mais groupe  fondamental du compl\'ementaire du n\oe ud de tr\`efle%
\footnote{Artin (\cite{A}, p. 53), voir aussi \cite{FLP} le corollaire p. 25
de l'expos\'e 2 de V. Po\'enaru.
}\
et, rappelant l'argument de Quillen pour construire un diff\'eomorphisme de
$M=SL(2, {\Bbb R})/SL(2, {\Bbb Z})$
sur le compl\'ementaire du n\oe ud de tr\`efle%
\footnote{L'espace des r\'eseaux de
${\bf R}^{2}\simeq{\bf C}$
s'identifie d'une part \`a
$GL(2, {\bf R})/GL(2, {\bf Z})$
et d'autre part, via l'\'equation 
$(W_{\Lambda})\!:(\wp')^{2}\!=\!4(\wp)^{3} - u\wp -v$
de la fonction de Weierstra\stz\ du r\'eseau
$\Lambda\subset{\bf C},$%
\hfill\break
$\wp(z)\!=\!{{1}\over{z^{2}}}+
\sum_{\lambda\in\Lambda\setminus\{0\}}{{1}\over{(z-\lambda)^{2}}}-{{1}\over{\lambda^{2}}},\,
u\!=\!g_{2}=60\sum_{\lambda\in\Lambda\setminus\{0\}}{{1}\over{\lambda^{4},}},\,
v\!=\!g_{3}=140\sum_{\lambda\in\Lambda\setminus\{0\}}{{1}\over{\lambda^{6}}}$%
, au compl\'ementaire dans
${\bf C}^{2}$
du discriminant
$\{\Delta_{\Lambda}(u, v)=27u^{2}-v^{3}=0\}$%
.
Ainsi l'espace des r\'eseaux unimodulaires
$SL(2, {\bf R})/SL(2, {\bf Z})$%
, puisque
$\Delta_{\Lambda}$
est quasi homog\`ene, s'identifie au compl\'ementaire
$\bigl({\bf C}^{2}\setminus\{\Delta_{\Lambda}=0\}\bigr)\cap\{|u|^{2}+|v|^{2}=1\}$
du n\oe ud de tr\`efle dans la sph\`ere unit\'e de
${\bf C}^{2}$%
.
},
obtient que le groupe fondamental de
$M$
est
$St(2, {\Bbb Z})\!\simeq\!B_{3}\!=\!\widetilde{SL}(2, {\Bbb Z})$%
,  pr\'eimage de
$SL(2, {\Bbb Z})$
dans le rev\^etement universel
$\widetilde{SL}(2, {\Bbb R})$
du groupe sp\'ecial lin\'eaire r\'eel
$SL(2, {\Bbb R})$%
, \'etablissant \`a nouveau la pr\'esentation tress\'ee d'Artin,
mais avec la r\'eserve
{\it (which is of course classical)}.
Sans doute pensait-il \`a un traitement
comme dans l'Appendice A de la monographie de C. Kassel et V. Turaev
(\cite{KT})
inspir\'e de celui de Reidemeister (\cite{Re})
de la pr\'esentation du groupe modulaire.

Voir
$B_{3}$
comme pr\'eimage de
$SL(2, {\Bbb Z})$
dans le rev\^etement universel de
$SL(2, {\Bbb R})$
faisait partie du folklore (voir Serre dans son cours de 1968-69 (\cite{Se}, p. 20).

L'int\'er\^et pour le th\'eor\`eme de Nielsen
${\romannumeral 3)\/}$
a \'et\'e raviv\'e depuis les ann\'ees 50 puis 80-90
par l'intervention du tore modulaire en arithm\'etique,
topologie des surfaces et 
d\'ecompte des g\'eod\'esiques simples
sur les surfaces hyperboliques%
\footnote{Voir \cite{Co}, l'article synth\`ese {\bf [H]} de Haas et la note \cite{McR}.
}.

S'ils attribuent le r\'esultat \`a Nielsen et
son article (\cite{N} de 1917), tant Magnus, Karass et Solitar (\cite{MKS}, p. 169)
que Lyndon et Shupp (\cite{LS}, p. 25),
 montrent%
 \footnote{Via r\'ef\'erence \`a Chang (\cite{Ch} pour Lyndon et Shupp, qui proposent aussi une
d\'emonstration fautive~: si
$Int(F_{2})\!=:\!I<\!K\!:=\!\ker(\rho)$%
, utilisant la pr\'esentation classique de
$GL(2, {\bf Z})$
avec un g\'e\-n\'e\-rateur d'ordre 2, un d'ordre 4 et un d'ordre 6 (%
$C=\begin{pmatrix}0& 1\cr 1&0\cr\end{pmatrix}$
et, avec nos notations,
$S, (AB)$%
) laissent le lecteur v\'erifier que les rel\`evements
$(v, u), (v^{-1}, u), (v^{-1}, uv)$
dans
$Aut(F_{2})$
de
$C, S, (AB)$
v\'erifient modulo
$I$
 les relations. M\^eme en corrigeant
$(v^{-1}, uv)$
en
$(v, u^{-1}v)$%
, cela montre que l'extension
$1\!\!\rightarrow\!\!K/I\!\!\rightarrow\!\!Aut(F_{2})/I\!\!\rightarrow\!\!GL(2, {\bf Z})\!\!%
\rightarrow\!\!1$
a une section, mais pas que
$K/I\!=\!\{1\}$%
!}\
 $ker(Aut(F_{2})\!\rightarrow\!GL(2, {\Bbb Z\/}))\!=\!Int(F_{2})$
en utilisant la pr\'esentation de
$Aut(F_{n})$ 
publi\'ee par Nielsen%
\footnote{{\sl Die Isomorphismengruppe der freien Gruppen\/} Math. Ann. {\bf 91\/} (1924), p. 169-209,\hfill\break
(traduction en anglais par J. Stillwell {\bf [N-12]}).
}\
sept ans apr\`es.

Le court (13 pages) article originel de Nielsen (\cite{N})
n'est pas de lecture ais\'ee~:

Nielsen signale que Dehn lui a communiqu\'e
${\romannumeral 4)\/}$
avec une preuve diff\'erente%
\footnote{sans r\'ef\'erence par Nielsen, imaginons~: la conservation 
\`a conjugaison
pr\`es du commutateur des g\'en\'erateurs par l'endomorphisme
$\varphi$
permet d'\'etendre en une
application de degr\'e
$1$
du tore
$f : {\bf R\/}^{2}/{\bf Z\/}^{2}=: R\rightarrow R$
sur lui-m\^eme  une application
$h\!: B\!\rightarrow\!B$
du bouquet
$B\!=\!\mu\cup\lambda\!:={\bf R\/}\!\times\!\{0\}/{\bf Z\/}\!\times\!\{0\}\cup%
\{0\}\!\times\!{\bf R\/}/\{0\}\!\times\!{\bf Z\/}$
sur lui-m\^ eme, identit\'e pr\`es de
$m\!=\!\overline{({{1}\over{2}}, {{1}\over{2}})}$
et induisant sur le groupe fondamental de
$B$
l'endomorphisme
$\varphi$%
.
On d\'eforme
$f:(R, R\!\setminus\!\{m\}, \{m\})\rightarrow (R, R\!\setminus\!\{m\}, \{m\})$
en une appliction lisse transverse \`a
$B$
et simplifie
$f^{-1}(\mu)$%
, puis
$f^{-1}(\lambda)$
(pour
$f^{-1}(\mu)$%
: on \'elimine  les composantes inessentielles puis, celles qui restent
(il y en a car
$f$
est de degr\'e
$1$%
)
sont, s'il y en a plus d'une,
des courbes parall\`eles, chaque paire de composantes voisines
\'etant s\'epar\'ee
par deux anneaux
$A_{1}, A_{2}$%
. Alors la restriction de
$f$
\`a
$A_{1}$%
, celui qui ne contient pas
$m$%
, est d'image disjointe de
$m$
donc de degr\'e nul et
$f$
se d\'eforme, pr\`es de
$A_{1}$%
, en une application envoyant
$A_{1}$
hors de
$\mu$%
\dots).
Ainsi
$f$
se d\'eforme en
$g$
lisse transverse \`a
$B$
avec
$g^{-1}(\mu)$
et
$g^{-1}(\lambda)$
connexes, donc
$g_{|}\!: g^{-1}(B)\!\rightarrow\!B$
est un rev\^etement connexe, de degr\'e
$deg(g_{|})=deg(g)=deg(f)=1$%
, ainsi 
$g_{|}$
est un hom\'eomorphisme, donc l'endomorphisme
$\varphi$
qu'il induit sur le groupe fondamental de
$\pi_{1}(R\!\setminus\!\{m\})$
est un isomorphisme.\hfill\findem
}.

Si l'analogue de
${\romannumeral 3)\/}$
et
${\romannumeral 4)\/}$
est clairement d\'egag\'e en fin de I et d\'ebut de II,
ce n'est que page 9  apr\`es une intrication de d\'emonstrations
et \'enonc\'es interm\'ediaires.

Le point
${\romannumeral 1)\/}$
est plut\^ot suppos\'e que d\'emontr\'e%
\footnote{\`a moins de savoir d\`es 1917 qu'un sous-groupe d'un groupe libre est libre
(attribu\'e \`a  Schreier (1927) ou Nielsen (1921) en Danois, traduction en anglais par A.W. et W.D.Neumann {\bf [N-8]}),
on ne peut affirmer avec Nielsen que si les g\'en\'erateurs
$a, b$
de
$F_{2}$
s'expriment en fonction de mots
$\alpha, \beta$
alors l'endomorphisme  envoyant
$a, b$ sur
$\alpha, \beta$
est un isomorphisme (pour avoir l'inverse, il faudrait  savoir que
$\alpha$
et
$\beta$
n'aient pas de relation entre eux).
}\
au tout d\'ebut de I.

N'ayant d\'egag\'e le sous-groupe
$H\!\subset\!Aut^{+}(F_{2})$%
, pour  l'analogue de
${\romannumeral 1)\/}$
et
${\romannumeral 3)\/}$%
, Nielsen proc\`ede en deux temps~:
d'abord de I 1. \`a I 6.  il r\'eduit jusqu'\`a 1,
par auto\-mor\-phisme int\'erieur et
 changement de \og paires images\fg
$(S(t_{1}), S(t_{2}))$
d\'ecrits dans le compl\'ement de l'{\petcap Affirmation 1\/} de l'Appendice B,
la longueur des images, par l'endomorphisme surjectif
$S$%
, des g\'en\'erateurs, 
puis de I 7. \`a I 10.,  examinant pas \`a pas le passage inverse et,
utilisant la derni\`ere partie de sa th\`ese%
\footnote{Kurvennetze au Fl\H achen Diss. Kiel 1913, p. 47-49 (accessible depuis 1986~:
{\bf [N-1]}).
}%
, il \'etablit que dans chacun des mots
$S(t_{i}),\, i\!=\!1, 2$
les exposants de chaque g\'en\'erateur ont m\^eme signe et que ces mots
ont une forme uniquement d\'etermin\'ee par l'ab\'elianis\'e de
$S$%
.

Cette derni\`ere partie a \'et\'e oubli\'ee%
\footnote{ou jug\'ee, sans acc\`es (avant 1986)
\`a la th\`ese de Nielsen, incompl\`ete.
}\
et a
provoqu\'e une grande activit\'e pour d\'eterminer la forme explicite des mots
$S(t_{i})$%
, avec McShane et Rivin, citons \cite{OZ}.

A la fin%
\footnote{\S9, p. 70 avec 
$\sigma_{1}, \sigma_{2}, a\!=\!\sigma_{1}\sigma_{2}, b\!=\!\sigma_{1}\sigma_{2}\sigma_{1}$
pour nos
$a, b, ab\!=:\!t, s$
 mais, comme nous
 $c\!=\!s^{2}$%
 !
}\
de \cite{A} Artin transforme sa pr\'esentation 
$B_{3}\!=< a, b\, |\, aba\!=\!bab>$
en
$<t, s\, |\, t^{3}\!=\!s^{2}> (t\!=\!ba, s\!=\!bab)$
et,  renvoyant \`a Dehn \cite{D} et Schreier \cite{Sc},
il met tout \'el\'ement
$w\!\in\!B_{3}$
sous la forme normale
$w\!=\!t^{\epsilon_{0}}st_{1}^{\epsilon_{1}}\cdots st^{\epsilon_{n-1}}s^{\epsilon_{n}}c^{k}$
avec~:
$$
n\in{\Bbb N\/}, k\!\in\!{\Bbb Z\/}, \epsilon_{j}\in\{-1,  1\}\
\hbox{\rm pour\/}\ 0< j < n,\, \epsilon_{0}\in\{-1,  0, 1\},\, \epsilon_{n}\in \{0, 1\}\ .
$$
Comme
$ts=babab=ca^{-1}, t^{-1}s=(ba)^{-1}bab=b$%
, quitte \`a modifier le traitement des termes ext\'erieurs~:
$w=s^{\nu_{1}}t^{\epsilon_{1}}s\cdots t^{\epsilon_{n-1}}ss^{\nu_{2}}c^{k},
\nu_{1}, \nu_{2}\!\in\!\{0, 1\}$%
, c'est celle  de l'Appendice C, obtenue par m\'ethode de Gau\stz\ sur les matrices
$2\times 2$
\`a c\oe fficients entiers naturels que l'on trouve \`a la fin de la th\`ese
de Nielsen {\bf [N-1]}.

La d\'etermination de la torsion de
$Aut(F_{2})$
est due \`a Meskin \cite{Me}. Il laisse cependant \`a la lectrice
une grande initiative dans la distinction entre mot r\'eduit de l'image par
un endomorphisme d'un mot r\'eduit et image de ce mot r\'eduit et
lui demande des
v\'erifications de ``non cancellation'' probl\'ematiques%
\footnote{par exemple p. 497 {\sl If
$W$
ends in
$B$%
, that is,
$W=UB$%
, then
$W^{S}=U^{S}A^{-1}B^{-1}$
and
$W^{S^{2}}=U^{S^{2}}A$
and the
$A'$
 do not cancel\/}
or, sans avoir prouv\'e pr\'ealablement que
$W$ ne peut commencer et finir par
$B$%
, si
$U=BU'$
commence aussi par
$B$%
, d'apr\`es la relation
$B^{S}=A^{-1}B^{-1}$%
, donn\'ee plus haut, il y a une simplification de
$A$
dans l'expression~:
$$(SW)^{3}=W^{S^{2}}W^{S}W (= (U^{S^{2}}A)(A^{-1}B^{-1}U'^{S})W)$$
}.
\vfill\eject


\def\refname{{\sc R\'ef\'erences}}


\vfill\eject



\section*{Quatri\`eme de couverture de la plaquette, sc\'enario de~:}

\section*{$SL(2, {\Bbb Z})$, les tresses \`a trois brins, le tore modulaire et
$Aut^{+}(F_{2})$}

The main subject of this booklet, written as it is
in french, is to propose an
``Artin presentation'' of $SL(2, {\Bbb Z})$.
This presentation is obtained via the action of the two-dimensional integer special linear group
on the so-called pillow orbifold ie the 2 sphere 
together with a singular euclidean metric,
having three cone points permuted by the action whilst 
the fourth cone point is fixed.

This orbifold is obtained as the quotient of the integer two-dimensional torus by
the elliptic involution $(x,y) \mapsto (-x,-y)$.
Artin's presentation of
$B_{3}$, braid group on three strings viewed as modular group of the sphere with
$3\!+\!1$
marked points, provides us with the desired presentation.

This Artin presentation provides a convenient description of the derived group
and its action on the Poincar\'e half plane,  an action whose quotient is
none other than the modular torus.

It appears that a formulation
of Nielsen's theorem for direct automorphisms of the free group
$F_{2}$
using Artin's braid group $B_{3}$
allows one to obtain in a simple and more rapid fashion the 
Artin presentation of
$SL(2, {\Bbb Z})$
which though surely known to experts
has not been popularized as it deserves.

However the proof of this formulation of Nielsen's theorem,
which seems not to be present in the literature at least
explicitly, needs Artin's presentation.
This is our task in the second appendix (an appendix which is 
not the core of the paper but merely intended as a kind of  flashback).

In  the remaining appendices we treat various subjects related to
Nielsen's theorem such as determining torsion elements of
$SL(2, {\Bbb Z\/})$
and the automorphism group of the free group of rank 2,
a corresponding formulation of Nielsen's theorem for the full automorphism group and the treatment of semi direct product of groups
with amalgamation  used in the formulation of Nielsen's theorem.

During the shooting of this short movie,
the ``director and script'', 
not only cited, but (re)read the
literature both old and new.
He hopes that this humble effort might aid,
the public  to avoid some of the many traps 
and to enjoy what he intends as
``a safari in the past and present ''.

\vfill

\noindent{
Metteur en sc\`ene et secr\'etaire~: Alexis Marin Bozonat\hfill\break
\null\ courriel~: alexis.charles.marin@gmail.com}

\noindent{
Institut Fourier, UMR 5582, Laboratoire de Math\'ematiques Universit\'e  Grenoble Alpes, CS 40700, 38058 Grenoble cedex 9, France}

\end{document}